\documentclass{elsarticle}

\usepackage{bm}
\usepackage{amsmath}
\usepackage{amsfonts}
\usepackage{amssymb}
\usepackage{graphicx}
\usepackage{paralist}
\usepackage{algorithm}
\usepackage{pstool}
\usepackage{wrapfig}
\usepackage{subcaption}
\usepackage{algpseudocode}
\usepackage{mathtools}
\usepackage{fullpage}

\usepackage{tikz}
\usepackage{pgfplots}
\usepackage{pgfplotstable, filecontents}
\pgfplotsset{compat=1.9}
\usetikzlibrary{pgfplots.groupplots}
\usepgfplotslibrary{fillbetween}
\usetikzlibrary{calc,fit,matrix,arrows,automata,positioning,shapes}
\usepackage{pgfplotstable, booktabs}
\usetikzlibrary{patterns,decorations.pathmorphing,decorations.markings,decorations.pathmorphing}
\pgfplotsset{select coords between index/.style 2 args={
    x filter/.code={
        \ifnum\coordindex<#1\fi
        \ifnum\coordindex>#2\fi
    }
}}

\pgfplotsset{compat=newest}

\usetikzlibrary{calc}

\newcommand{\pder}[2]{\ensuremath{\frac{\partial #1}{\partial #2}}} 

\newcommand{\oder}[2]{\ensuremath{\frac{\mathrm{d} #1}{\mathrm{d} #2}}} 

\newcommand{\Fcal}{\ensuremath{\mathcal{F}}}
\newcommand{\Gcal}{\ensuremath{\mathcal{G}}}

\newcommand{\Jcal}{\ensuremath{\mathcal{J}}}

\newcommand{\Lcal}{\ensuremath{\mathcal{L}}}

\newcommand\Mbm{{\ensuremath{\bm{M}}}}

\newcommand\Rbm{{\ensuremath{\bm{R}}}}

\newcommand\cbm{{\ensuremath{\bm{c}}}}

\newcommand\fbm{{\ensuremath{\bm{f}}}}
\newcommand\gbm{{\ensuremath{\bm{g}}}}

\newcommand\kbm{{\ensuremath{\bm{k}}}}

\newcommand\pbm{{\ensuremath{\bm{p}}}}
\newcommand\qbm{{\ensuremath{\bm{q}}}}
\newcommand\rbm{{\ensuremath{\bm{r}}}}

\newcommand\ubm{{\ensuremath{\bm{u}}}}

\newcommand\lambdabold{{\ensuremath{\boldsymbol{\lambda}}}}

\newcommand\mubold{{\ensuremath{\boldsymbol{\mu}}}}
\newcommand\kappabold{{\ensuremath{\boldsymbol{\kappa}}}}

\newcommand\taubold{{\ensuremath{\boldsymbol{\tau}}}}

\newcommand\sigmabold{{\ensuremath{\boldsymbol{\sigma}}}}

\newcommand\zerobold{\ensuremath{\mathbf{0}}}

\usepackage{xifthen}
\usepackage[thinlines]{easytable}
\newboolean{fastcompile}
\setboolean{fastcompile}{false}

\newcommand{\dmu}[1][]{\ifthenelse{\isempty{#1}}{\partial_\mubold}{\pder{#1}{\mubold}}}
\newcommand{\mass}[1][]{\ifthenelse{\isempty{#1}}{\Mbm}{\Mbm^{#1}}}
\newcommand{\stvc}[1][]{\ifthenelse{\isempty{#1}}{\ubm}{\ubm^{#1}}}
\newcommand{\stvcdot}[1][]{\ifthenelse{\isempty{#1}}{\dot\ubm}{\dot\ubm^{#1}}}
\newcommand{\stvcbar}[1][]{\ifthenelse{\isempty{#1}}{\bar\ubm}{\bar\ubm^{#1}}}
\newcommand{\res}[1][]{\ifthenelse{\isempty{#1}}{\rbm}{\rbm^{#1}}}
\newcommand{\resimpl}[1][]{\ifthenelse{\isempty{#1}}{\gbm}{\gbm^{#1}}}
\newcommand{\resexpl}[1][]{\ifthenelse{\isempty{#1}}{\fbm}{\fbm^{#1}}}
\newcommand{\cpl}[1][]{\ifthenelse{\isempty{#1}}{\cbm}{\cbm^{#1}}}
\newcommand{\cplprd}[1][]{\ifthenelse{\isempty{#1}}{\tilde\cbm}{\tilde\cbm^{#1}}}

\newcommand{\pstp}[2][]{\ifthenelse{\isempty{#2}}{\ubm_{#1}}{\ubm_{#2}^{#1}}}
\newcommand{\pstgki}[3][]{\ifthenelse{\isempty{#3}}{\kbm_{#1,#2}}{\kbm_{#2,#3}^{#1}}}
\newcommand{\pstgke}[3][]{\ifthenelse{\isempty{#3}}{\hat\kbm_{#1,#2}}{\hat\kbm_{#2,#3}^{#1}}}
\newcommand{\pstgu}[3][]{\ifthenelse{\isempty{#3}}{\ubm_{#1,#2}}{\ubm_{#2,#3}^{#1}}}
\newcommand{\cplstp}[2][]{\ifthenelse{\isempty{#2}}{\cbm_{#1}}{\cbm_{#2}^{#1}}}
\newcommand{\cplprdstp}[2][]{\ifthenelse{\isempty{#2}}{\tilde\cbm_{#1}}{\tilde\cbm_{#2}^{#1}}}

\newcommand{\dt}[1]{\Delta t_{#1}}

\newcommand{\tstgi}[2]{t_{#1,#2}}

\newcommand{\rstp}[2][]{\ifthenelse{\isempty{#2}}{\tilde\rbm_{#1}}{\tilde\rbm_{#2}^{#1}}}  
\newcommand{\rstgi}[3][]{\ifthenelse{\isempty{#3}}{\Rbm_{#1,#2}}{\Rbm_{#2,#3}^{#1}}}
\newcommand{\rstge}[3][]{\ifthenelse{\isempty{#3}}{\hat\Rbm_{#1,#2}}{\hat\Rbm_{#2,#3}^{#1}}}
\newcommand{\qstg}[3][]{\ifthenelse{\isempty{#3}}{\qbm_{#1,#2}}{\qbm_{#2,#3}^{#1}}}
\newcommand{\prstp}[2][]{\ifthenelse{\isempty{#2}}{\pbm_{#1}}{\pbm_{#2}^{#1}}}

\newcommand{\dstp}[1]{\lambdabold_{#1}}
\newcommand{\dstgi}[2]{\kappabold_{#1,#2}}
\newcommand{\dstge}[2]{\hat\kappabold_{#1,#2}}
\newcommand{\sigstp}[1]{\sigmabold_{#1}}
\newcommand{\taustg}[2]{\taubold_{#1,#2}}

\newcommand{\eqnref}[1]{Eq. (\ref{#1})}
\newcommand{\figref}[1]{Figure \ref{#1}}
\newcommand{\tabref}[1]{Table \ref{#1}}

\usepackage{graphicx}
\usepackage[version=4]{mhchem}
\usepackage{siunitx}
\usepackage{longtable,tabularx}
\usepackage{booktabs}
\title{A high-order partitioned solver for general multiphysics problems and its applications in optimization}

\author[rvt1]{D.~Z.~Huang\fnref{fn1}\corref{cor1}}
\ead{zhengyuh@stanford.edu}

\author[rvt2,rvt3]{P.-O.~Persson\fnref{fn2}}
\ead{persson@berkeley.edu}

\author[rvt2,rvt4]{M.~J.~Zahr\fnref{fn3,fn4}}
\ead{mjzahr@lbl.gov}

\address[rvt1]{Institute for Computational and Mathematical Engineering,
               Stanford University, Stanford, CA, 94305, United States}
\address[rvt2]{Mathematics Group, Lawrence Berkeley National Laboratory,
               1 Cyclotron Road, Berkeley, CA 94720, United States}
\address[rvt3]{Department of Mathematics, University of California, Berkeley,
               Berkeley, CA 94720, United States}
\address[rvt4]{Department of Aerospace and Mechanical Engineering,
               University of Notre Dame, Notre Dame, IN, 46556, United States}
\cortext[cor1]{Corresponding author}

\fntext[fn1]{Graduate Student, Institute for Computational and Mathematical
             Engineering, Stanford University}
\fntext[fn2]{Associate Professor, Department of Mathematics, University of
             California, Berkeley}
\fntext[fn3]{Luis W. Alvarez Postdoctoral Fellow,
             Computational Research Division,
             Lawrence Berkeley National Laboratory}
\fntext[fn4]{Assistant Professor, Department of Aerospace and Mechanical
             Engineering, University of Notre Dame}

\begin{document}

\begin{abstract}
A high-order accurate adjoint-based optimization framework  is presented for
unsteady multiphysics problems.  The fully discrete adjoint solver relies
on the  high-order, linearly stable, partitioned solver introduced in
\cite{huang2018high},  where different subsystems are modeled and
discretized separately.  The coupled system of semi-discretized ordinary
differential equations is taken as a monolithic system and partitioned using
an implicit-explicit Runge-Kutta (IMEX-RK) discretization
\cite{ascher1997implicit}. Quantities of interest (QoI) that take the form
of space-time integrals are discretized in a solver-consistent manner. 
The corresponding adjoint equations are derived to compute exact gradients of
QoI,  which can be solved in a partitioned manner,  i.e. subsystem-by-subsystem
and substage-by-substage, thanks to the partitioned primal solver. 
These quantities of interest and their gradients are then used in the context
of gradient-based PDE-constrained optimization. The present optimization
framework is applied to two fluid-structure interaction problems: 1D piston
problem with a three-field formulation and a 2D energy harvesting problem with
a two-field formulation.
\end{abstract}
\maketitle
\section{Introduction}
Optimization problems involving multiphysics systems commonly arise in engineering practice, particularly in the context of design or 
control of physics-based systems. These problems lead to PDE-constrained optimizations.
In the literature, a majority of research in PDE-constrained optimization has been focused on a single physical system or steady PDEs,  
which is sufficient for a large class of problems of interest.
However, there is a large class of problems where such analysis is insufficient, such as problems that involves the interactions of multiple physical 
systems or physical phenomena, which are generally inherently dynamic. Typical examples include flapping flight for Micro-Aerial Vehicles (MAVs) designs \cite{chen2007numerical, stanford2010analytical}, 
optimal combustion control system to maintain stable combustion with low exhaust emissions \cite{day2000numerical, jeong2006optimization}, microscale swimmer designs 
for drug delivery \cite{dreyfus2005microscopic}, and wind turbine performance optimization \cite{mosetti1994optimization, xudong2009shape} to extract maximum energy. Design and control of these types of systems are
challenging considering the coupling effects of multiple physics and the high computational cost due to their unsteady nature. Innovative multiphysics solvers 
and state-of-art optimization tools are needed to solve such problems.

We first review the high-order, linearly stable, implicit-explicit Runge-Kutta (IMEX-RK) \cite{ascher1997implicit} based partitioned solvers 
for multiphysics problems proposed in \cite{huang2018high}.
In this framework, a generic multiphysics problem is modeled as a system of $n$ 
systems of partial differential equations where the $i$th subsystem is coupled to the other subsystems through a coupling term that can depend on
the state of all the other subsystems. This coupled system of partial differential equations reduces to a coupled system of ordinary differential equations
via the method of lines where an appropriate spatial discretization is applied to each subsystem. The coupled system of ordinary differential
equations is taken as a monolithic system and discretized using an IMEX-RK discretization with a specific implicit-explicit decomposition that
introduces the concept of a \emph{predictor} for the coupling term. 
Four coupling predictors are proposed that enable the monolithic system to be solved in a partitioned manner, i.e., subsystem-by-subsystem,
and preserve the IMEX-RK structure and therefore the design order of accuracy of the monolithic
scheme. The four partitioned solvers that result from these predictors are high-order accurate, allow for
maximum re-use of existing single-physics software, and two of the four solvers allow the subsystems to
be solved in parallel at a given stage and time step. In \cite{huang2018high}, we also analyze the stability of a coupled, linear
model problem and show that one of the partitioned solvers achieves unconditional linear stability, while
the others are unconditionally stable only for certain values of the coupling strength. 

Next, we derive the corresponding fully discrete sensitivity and adjoint equations for general optimization problems. 
Here, we mainly focus on one of the aforementioned partitioned solvers, the weakly coupled Gauss-Seidel predictor based partitioned solver, which has demonstrated 
its high-order accuracy, numerical stability and software maintainability in many engineering problems.
Quantities of interest or objective functions, e.g.  energy consumption or the quantities of combustion emission, 
that take the form of space-time integrals that are discretized in a solver-consistent manner.
This ensures the discretization order of quantities of interest exactly matches the PDE temporal discretization.
The aforementioned multiphysics partitioned solver becomes the PDE-constraint of the optimization problem. 
To compute exact gradients of quantities of interest, we need to solve the multiphysics problem, and then either 
compute the sensitivity of the state variables through forward time-marching or evaluate the 
adjoint variables through backward time-marching. We can leverage the high-order linear stability property of the partitioned solver, 
which takes large time steps and therefore reduces the number of time steps and accelerates the time-marching procedure. 
The optimization solver IPOPT \cite{wachter2006implementation} is 
used to solve the optimization problem based on a nonlinearly constrained interior point method.

The remainder of this paper is organized as follows. In Section~\ref{SEC: EQUATIONS}, the governing equations of the multiphysics system, the integral form quantities of interest and their semi-discretizations are introduced. The high-order temporal discretization, based on IMEX-RK schemes, is described in Section~\ref{SEC: PARTITIONED SOLVER}, which leads to a partitioned multiphysics solver.  Following this in Section~\ref{SEC: OPTIMIZATION}, the corresponding fully discrete sensitivity equations and adjoint equations are derived, which deliver the \emph{exact} gradient of the QoIs. Section~\ref{SEC: APPLICATIONS} demonstrates the approach as applied to two optimization problems: a 1D oscillating piston problem and a 2D airfoil energy harvesting problem. Finally, conclusions are offered in Section~\ref{SEC: CONCLUSION}.

\section{Governing multiphysics equations and semi-discretization}\label{SEC: EQUATIONS}
Consider a general formulation of a mathematical model describing the
behavior of multiple interacting physical phenomena described by the
following coupled system of partial differential equations
\begin{equation} \label{EQ: GOVERN}
 \partial_t u^i = \Lcal^i(u^i,\,c^i,\,x,\,\mubold,\,t),
 \quad x\in\Omega^i(c^i,\,\mubold,\,t), \quad t \in (0,\,T)
\end{equation}
for $i = 1,\,\dots,\,m$, where $m$ represents the number of physical systems,
and boundary conditions are excluded for brevity.
The $i$th physical system is modeled as a partial differential equation
characterized by the generalized differential operator $\Lcal^i$ that
defines a conservation law or other type of balance law, the state variable
$u^i$ that is the solution of the $i$th physical system on the
space-time domain $\Omega^i \times (0,\,T)$, and a \emph{coupling term} $c^i$
that, in general, couples the $i$th system to the other $m-1$ systems. 
In the general case, the differential operator $\Lcal^i$, domain $\Omega^i$,
and boundary conditions depend on the coupling term.
The coupling term contains quantities usually considered \emph{data} required
to define the $i$th PDE, such as boundary conditions or material properties.
In a single-physics setting, these quantities would be prescribed, but in the
multiphysics setting they are determined from the state vectors of all $m$
systems, i.e.,
\begin{equation} \label{EQ: CPLTERM}
 c^i = c^i(u^1,\,\dots,\,u^m,\,x,\,\mubold,\,t).
\end{equation}
The definition of the coupling term is problem-dependent and special structure 
in the coupling term can be exploited to create a
better partitioned solver. While the form of (\ref{EQ: GOVERN}) is specific
to first-order temporal systems, it includes equations with higher-order
temporal derivatives, assuming they have been re-cast in the first-order form.
The spatial domains $\Omega^i$ for the individual systems may or may not be
overlapping and in many cases are the same, i.e., $\Omega^i = \Omega$
for $i = 1,\,\dots,\,m$. 
Any of the operators, solution variables, and even the deformed 
computational domain might depend on the parameter vector $\mubold$.
The quantities of interest are assumed to be of the integral form, 
\begin{equation}
\label{EQ: QoI}
 \Jcal(u^1,\,\dots,\,u^m,\,x,\,\mubold,\,T) = \int_0^{T}j(u^1(\tau),\,\dots,\,u^m(\tau),\,\mubold,\,\tau)d\tau.
\end{equation}

We introduce the semi-discrete
form of the coupled partial differential equations in (\ref{EQ: GOVERN}) that
arises from applying an appropriate spatial discretization to the $i$th PDE
system individually, which takes the form
\begin{equation} \label{EQ: GOVERN-SD-INDIV}
 \mass[i] \stvcdot[i] = \res[i](\stvc[i],\,\cpl[i],\,\mubold,\,t),
 \quad t \in (0,\,T)
\end{equation}
where $\stvc[i](t)$ is the semi-discrete state vector corresponding to the
spatial discretization of $u^i(x,\,t)$, $\res^i$ is the spatial
discretization of the differential operator $\Lcal^i$ and called the
velocity of the ODE system in the remainder of the document, and $\cpl[i]$ is
the semi-discrete coupling term corresponding to the spatial discretization of
$c^i(u^1,\,\dots,\,u^m,\,x,\,t)$. In general, the coupling term depends on the
semi-discrete state vector of all $m$ systems
\begin{equation} \label{EQ: CPLTERM-SD-INDIV}
 \cpl[i] = \cpl[i](\stvc[1],\,\dots,\,\stvc[m],\,\mubold,\,t).
\end{equation}
For convenience, we re-write the system of ordinary differential equations
in (\ref{EQ: GOVERN-SD-INDIV})-(\ref{EQ: CPLTERM-SD-INDIV}) as
\begin{equation} \label{EQ: GOVERN-SD}
 \mass\stvcdot = \res(\stvc,\,\cpl(\stvc,\,\mubold,\,t),\,\mubold,\,t),
 \quad t \in (0,\,T),
\end{equation}
where the combined mass matrix is a block diagonal matrix consisting of
the single-physics mass matrices
\begin{equation*}
 \mass = \begin{bmatrix} \mass[1] && \\ & \ddots & \\ && \mass[m] \end{bmatrix}
\end{equation*}
and the combined state vector, coupling term, and nonlinear residual are
vectors consisting of the corresponding single-physics term, concatenated
across all $m$ systems
\begin{equation*}
 \stvc = \begin{bmatrix} \stvc^1 \\ \vdots \\ \stvc^m \end{bmatrix} \qquad
 \cpl(\stvc,\,t) = \begin{bmatrix}
                     \cpl^1(\stvc^1,\,\dots,\,\stvc^m,\,\mubold,\,t) \\
                     \vdots \\
                     \cpl^m(\stvc^1,\,\dots,\,\stvc^m,\,\mubold,\,t)
                  \end{bmatrix} \qquad
 \res(\stvc,\,\cpl,\,t) = \begin{bmatrix}
                            \res^1(\stvc^1,\,\cpl^1,\,\mubold,\,t) \\
                            \vdots \\
                            \res^m(\stvc^m,\,\cpl^m,\,\mubold,\,t)
                          \end{bmatrix}.
\end{equation*}
The total derivative, or Jacobian, of the semi-discrete velocity
$D_{\stvc}\res$ is expanded as
\begin{equation}\label{EQ: GOVERN-Jac}
 D_{\stvc}\res = \pder{\res}{\stvc} + \pder{\res}{\cpl}\pder{\cpl}{\stvc},
\end{equation}
where the individual terms take the form
\begin{equation} \label{EQ: drdu drdc dcdu}
 \pder{\res}{\stvc} =
 \begin{bmatrix} \displaystyle{\pder{\res[1]}{\stvc[1]}} & & \\
                 & \ddots & \\
                 & & \displaystyle{\pder{\res[m]}{\stvc[m]}}
 \end{bmatrix} \qquad
 \pder{\res}{\cpl} =
 \begin{bmatrix} \displaystyle{\pder{\res[1]}{\cpl[1]}} & & \\
                 & \ddots & \\
                 & & \displaystyle{\pder{\res[m]}{\cpl[m]}}
 \end{bmatrix} \qquad
 \pder{\cpl}{\stvc} =
 \begin{bmatrix} \displaystyle{\pder{\cpl[1]}{\stvc[1]}} & \cdots &
                 \displaystyle{\pder{\cpl[1]}{\stvc[m]}} \\
                 \vdots & \ddots & \vdots \\
                 \displaystyle{\pder{\cpl[m]}{\stvc[1]}} & \cdots &
                 \displaystyle{\pder{\cpl[m]}{\stvc[m]}}
 \end{bmatrix},
\end{equation}
and the dependencies have been dropped for brevity.
The first term in the Jacobian, \eqnref{EQ: GOVERN-Jac}, is block diagonal and accounts for the direct
contribution of a state to its own system while the second term accounts
for the coupling between systems. 

\section{A high-order partitioned solver for multiphysics problems}
\label{SEC: PARTITIONED SOLVER}
In this section, high-order partitioned time-integration schemes for multiphysics systems are introduced. 
In a partitioned sense,  individual off-the-shelf single-physics solvers are combined to solve the 
multiphysics problem, rather than considering the monolithic multiphysics system. However, 
they tend to be limited to low-order accuracy and
have stringent stability requirements. Our partitioned time-integration
scheme mitigates most of these issues by combining high-order
implicit-explicit Runge-Kutta (IMEX) schemes for the monolithic multiphysics
system with a judicious implicit-explicit decomposition that \emph{diagonally}
couples the individual systems via a novel predictor for the coupling terms.

\subsection{Background: implicit-explicit Runge-Kutta schemes}
Implicit-explicit Runge-Kutta schemes, first proposed in \cite{zhong1996additive, ascher1997implicit}, define a
family of high-order discretizations for nonlinear differential equations
whose velocity term can be decomposed into a sum of a non-stiff $\resexpl$
and stiff $\resimpl$ velocity
\begin{equation} \label{EQ: ODEADDSPLIT}
 \mass\dot\stvc = \resexpl(\stvc,\,t) + \resimpl(\stvc,\,t).
\end{equation}
The non-stiff $\resexpl$ velocity is integrated with an $s$-stage explicit
Runge-Kutta scheme and the stiff term $\resimpl$ is integrated with an
$s$-stage diagonally implicit Runge-Kutta scheme.
IMEX Runge-Kutta schemes are compactly represented by a double tableau in
the usual Butcher notation (\tabref{TAB: IMEX}), where $\hat{A}$, $\hat{b}$,
$\hat{c}$ defines the Butcher tableau for the explicit Runge-Kutta
scheme used for $\resexpl$ and $A$, $b$, $c$ defines the diagonally implicit
Runge-Kutta scheme used for $\resimpl$.
In this work, we mainly consider IMEX-RK schemes proposed in \cite{christopher2001additive}, in which the
implicit Runge-Kutta part of these IMEX schemes are L-stable,
stiffly-accurate, and have an explicit first stage ($a_{11} = 0$).

\begin{table}[!htb]
\begin{minipage}{.5\linewidth}
\caption*{Explicit Runge-Kutta coefficients}
\centering
\begin{tabular}{c|ccccc}
0\\
$\hat c_2$    & $\hat a_{21}$\\
$\hat c_3$    & $\hat a_{31}$ & $\hat a_{32}$\\
$\vdots$ & $\vdots$ &           & $\ddots$\\
$\hat c_s$    & $\hat  a_{s1}$ & $a_{s2}$  & $\cdots$ & $\hat a_{ss - 1}$\\
\hline
& $\hat b_1$    & $\hat b_2$     & $\cdots$ & $\hat b_{s-1}$     & $\hat b_s$\\
\end{tabular}
\end{minipage}
\begin{minipage}{.5\linewidth}
\centering
\caption*{Implicit Runge-Kutta coefficients}
 \begin{tabular}{c|ccccc}
$c_1$\\
$c_2$    & $a_{21}$ & $a_{22}$\\
$c_3$    & $a_{31}$ & $a_{32}$ & $a_{33}$\\
$\vdots$ & $\vdots$ &      &     & $\ddots$\\
$c_s$    & $a_{s1}$ & $a_{s2}$  & $\cdots$ & $a_{ss-1}$ & $a_{ss}$\\
\hline
& $b_1$    & $b_2$     & $\cdots$ & $b_{s-1}$     & $b_s$\\
\end{tabular}
 \end{minipage} 
\caption{Butcher Tableaux for an $s$-stage implicit-explicit Runge-Kutta scheme.}   
\label{TAB: IMEX} 
\end{table}

Consider a discretization of the time domain $[0,\,T]$ into $N_t$ segments
with endpoints $\{t_0,\,\dots,\,t_{N_t}\}$, with the $n$th segment having
length $\Delta t_n = t_n - t_{n-1}$ for $n = 1,\,\dots,\,N_t$. Also,
let $\pstp{n}$ denote the approximation of the solution of the
differential equation in (\ref{EQ: ODEADDSPLIT}) at timestep $n$, i.e.,
$\pstp{n} \approx \stvc(t_n)$. Then, given the explicit
$(\hat{A},\,\hat{b},\,\hat{c})$ and implicit $(A,\,b,\,c)$ Butcher tableaux,
the $s$-stage IMEX Runge-Kutta scheme that advances $\pstp{n-1}$ to
$\pstp{n}$ is given by
\begin{subequations}
\label{EQ: IMEX STEP IMP EXP STAGE}
 \begin{align}
  \pstp{n} &= \pstp{n-1} + \sum_{p=1}^s\hat{b}_p\pstgke{n}{p} +
                            \sum_{p=1}^s b_p\pstgki{n}{p},
                                                        \label{EQ: IMEX STEP} \\
  \mass\pstgki{n}{j} &= \dt{n}\resimpl(\pstgu{n}{j},\,t_{n-1}+c_j\dt{n}),
                                                         \label{EQ: IMEX IMP} \\
  \mass\pstgke{n}{j} &= \dt{n}\resexpl(\pstgu{n}{j},\,t_{n-1}+\hat{c}_j\dt{n}),
                                                         \label{EQ: IMEX EXP} \\
  \pstgu{n}{j} &= \pstp{n-1} + \sum_{p=1}^{j-1}\hat{a}_{jp}\pstgke{n}{p} +
                               \sum_{p=1}^j a_{jp}\pstgki{n}{p},
                                                         \label{EQ: IMEX STAGE}
 \end{align}
\end{subequations}
where $\pstgke{n}{p}$ and $\pstgki{n}{p}$ are the $p$th explicit and implicit
velocity stage, respectively, corresponding to timestep $n$ and
$\pstgu{n}{p}$ is the approximation to $\pstp{n}$ at stage $p$ of timestep $n$. For each
stage $j$, the nonlinear system of equations in (\ref{EQ: IMEX IMP}) must
be solved to compute the implicit stage $\pstgki{n}{j}$. Next, the explicit
stage can be computed directly from (\ref{EQ: IMEX EXP}) since the stage
approximation $\pstgu{n}{j}$ does not depend on the explicit stage
$\pstgke{n}{j}$. Finally, given the previous timestep and all implicit
and explicit stages, the solution at time $n$ is determined from
(\ref{EQ: IMEX STEP}).

\subsection{A partitioned implicit-explicit Runge-Kutta scheme for
            multiphysics systems}
The proposed high-order partitioned scheme for integration of generic
time-dependent multiphysics problems of the form
(\ref{EQ: GOVERN-SD-INDIV})-(\ref{EQ: CPLTERM-SD-INDIV}) is built on an
IMEX Runge-Kutta discretization of the monolithic system.
A special choice of implicit-explicit decomposition, along with the
introduction of predictors for the coupling term, creates a
\emph{diagonal} or \emph{triangular} dependency between the
systems and allows the monolithic discretization to be solved in
a partitioned manner. 
The proposed decomposition handles a majority of
the relevant physics \emph{implicitly} to leverage the enhanced stability
properties of such schemes, while only the correction to the coupling
predictor is handled explicitly. 

\subsubsection{Implicit-explicit decomposition and monolithic IMEX
               Runge-Kutta discretization}
To begin our construction, recall the semi-discrete form of the multiphysics
system (\ref{EQ: GOVERN-SD}) and consider the splitting of the
velocity term $\res(\stvc,\,\cpl(\stvc,\,t),\,t)$ as
\begin{equation} \label{EQ: RES SPLIT}
 \res(\stvc,\,\cpl(\stvc,\,\mubold,\,t),\,\mubold,\,t) = \resexpl(\stvc,\,\cplprd,\,\mubold,\,t) +
                                     \resimpl(\stvc,\,\cplprd,\,\mubold,\,t)
\end{equation}
where $\cplprd$ is an approximation, or \emph{predictor}, of the coupling
term $\cpl(\stvc,\,t)$ and the terms that will be handled explicitly
$\resexpl$ and implicitly $\resimpl$ in the IMEX discretization are defined as
\begin{subequations} \label{EQ: IMEX MP SPLIT}
 \begin{align}
  \resexpl(\stvc,\,\cplprd,\,\mubold,\,t) &= \res(\stvc,\,\cpl(\stvc,\,\mubold,\,t),\,\mubold,\,t) - 
                                   \res(\stvc,\,\cplprd,\,\mubold,\,t)
                                   \label{EQ: IMEX MP EXP} \\
  \resimpl(\stvc,\,\cplprd,\,\mubold,\,t) &= \res(\stvc,\,\cplprd,\,\mubold,\,t),
                                   \label{EQ: IMEX MP IMP}
 \end{align}
\end{subequations}
where the dependence on the predictor is explicitly included. In general,
the predictor depends on the instantaneous state vector $\stvc(t)$ and
data $\stvcbar$, likely from the history of the state vector
$\{\stvc(\tau)\,|\,\tau < t\}$
\begin{equation} \label{EQ: CPLPRD GEN}
 \cplprd = \cplprd(\stvc,\,\stvcbar,\,\mubold,\,t).
\end{equation}

With this decomposition of the velocity of the semi-discrete multiphysics
system in (\ref{EQ: IMEX MP SPLIT}), the IMEX Runge-Kutta scheme in
(\ref{EQ: IMEX STEP IMP EXP STAGE}) applied to the monolithic multiphysics
system (\ref{EQ: GOVERN-SD}) becomes
\begin{subequations} 
 \begin{align}
  \pstp{n} &= \pstp{n-1} + \sum_{p=1}^s\hat{b}_p\pstgke{n}{p} +
                            \sum_{p=1}^s b_p\pstgki{n}{p}, \label{EQ: IMEX STEP IMP EXP STAGE WITH PRED 1}\\
  \mass\pstgki{n}{j} &= \dt{n}\resimpl(\pstgu{n}{j},\,
                                   \cplprd(\pstgu{n}{j},\,\pstp{n-1},\,
                                           \mubold,\,\tstgi{n}{j}),\,
                                   \mubold,\,\tstgi{n}{j}), \label{EQ: IMEX STEP IMP EXP STAGE WITH PRED 2}\\
  \mass\pstgke{n}{j} &= \dt{n}\resexpl(\pstgu{n}{j},\,
                                   \cplprd(\pstgu{n}{j},\,\pstp{n-1},\,
                                           \mubold,\,\tstgi{n}{j}),\,
                                   \mubold,\,\tstgi{n}{j}), \label{EQ: IMEX STEP IMP EXP STAGE WITH PRED 3}\\
  \pstgu{n}{j} &= \pstp{n-1} + \sum_{p=1}^{j-1}\hat{a}_{jp}\pstgke{n}{p} +
                               \sum_{p=1}^j a_{jp}\pstgki{n}{p}, \label{EQ: IMEX STEP IMP EXP STAGE WITH PRED 4}
 \end{align}
\end{subequations}
where the data used in the coupling predictor is taken from the previous
timestep. This is the general form of the fully discrete, monolithic
multiphysics system where the coupling predictor is unspecified.
In the general setting where each coupling predictor depends on the state
of all systems, the Jacobian of the coupling predictor is block
dense with potentially sparse blocks
\begin{equation*}
 \pder{\cplprd}{\stvc} =
 \begin{bmatrix} \displaystyle{\pder{\cplprd[1]}{\stvc[1]}} & \cdots &
                 \displaystyle{\pder{\cplprd[1]}{\stvc[m]}} \\
                 \vdots & \ddots & \vdots \\
                 \displaystyle{\pder{\cplprd[m]}{\stvc[1]}} & \cdots &
                 \displaystyle{\pder{\cplprd[m]}{\stvc[m]}}
 \end{bmatrix}.
\end{equation*}
This implies the Jacobian of the implicit velocity
\begin{equation*}
 D_{\stvc} \resimpl = \pder{\res}{\stvc} +
                      \pder{\res}{\cplprd} \pder{\cplprd}{\stvc}
\end{equation*}
is also block dense, which highlights the fact that there is coupling across
all systems and a monolithic solver is required for the implicit step.

\subsubsection{Weakly coupled Gauss-Seidel predictor}
The Gauss-Seidel-type (triangular) predictors for the multiphysics system
assume the individual systems are \emph{ordered} in a physically relevant
manner. The preferred ordering is problem-dependent.
The weakly coupled Gauss-Seidel-type predictor for the $i$th system is
defined as
\begin{equation} \label{EQ: CPLPRD WEAK GS}
 \cplprd[i](\stvc,\,\stvcbar,\,\mubold,\,t) = \cpl(\stvc[1],\,\dots,\,\stvc[i-1],\,
                                     \stvcbar[i],\,\dots,\,\stvcbar[m],\,\mubold,\,t)
\end{equation}
for $i = 1,\,\dots,\,m$. At the fully discrete level, this predictor takes
the form
\begin{equation} \label{EQ: CPLPRD WEAK GS FD}
 \cplprd[i](\pstgu{n}{j},\,\pstp{n-1},\,\mubold,\,t) =
                             \cpl(\pstgu[1]{n}{j},\,\dots,\,\pstgu[i-1]{n}{j},\,
                                  \pstp[i]{n-1},\,\dots,\,\pstp[m]{n-1},\,\mubold,\,t).
\end{equation}
In the context of the IMEX-RK discretization in
(\ref{EQ: IMEX STEP IMP EXP STAGE WITH PRED 1}-\ref{EQ: IMEX STEP IMP EXP STAGE WITH PRED 4}), the $i$th predictor
lags the state of systems $i,\,\dots,\,m$ to the previous timestep in the
evaluation of the coupling term throughout all stages of the timestep.
The IMEX-RK discretization of the multiphysics system in
(\ref{EQ: IMEX STEP IMP EXP STAGE WITH PRED 1}-\ref{EQ: IMEX STEP IMP EXP STAGE WITH PRED 4}) with this form of the predictor leads to  Algorithm~\ref{ALG: IMEX WEAK GS}.
\begin{algorithm}
 \caption{Implicit-Explicit Runge-Kutta partitioned multiphysics scheme:
          weak Gauss-Seidel predictor}
 \label{ALG: IMEX WEAK GS}
 \begin{algorithmic}[1]
   \For{stages $j = 1,\,\dots,\,s$}
    \For{physical systems $i = 1,\,\dots,\,m$}
     \State Define stage solution according to (\ref{EQ: IMEX STEP IMP EXP STAGE WITH PRED 1}):
      $\displaystyle{\pstgu[i]{n}{j} =
                                  \pstp[i]{n-1} +
                                  \sum_{p=1}^{j-1}\hat{a}_{jp}\pstgke[i]{n}{p} +
                                  \sum_{p=1}^j a_{jp}\pstgki[i]{n}{p}}$
     \State Implicit solve (\ref{EQ: IMEX STEP IMP EXP STAGE WITH PRED 2}) for $\pstgki[i]{n}{j}$:
      $\displaystyle{\mass[i]\pstgki[i]{n}{j} =
                     \dt{n}\resimpl[i](\pstgu[i]{n}{j},\,
                                    \cpl[i](\pstgu[1]{n}{j},\,\dots,\,
                                            \pstgu[i-1]{n}{j},\,
                                            \pstp[i]{n-1},\,\dots,\,
                                            \pstp[m]{n-1},\,
                                            \mubold,\,\tstgi{n}{j}),\,
                                    \mubold,\,\tstgi{n}{j})}$
     \State Explicit solve (\ref{EQ: IMEX STEP IMP EXP STAGE WITH PRED 3}) for $\pstgke[i]{n}{j}$:
      $\displaystyle{\mass[i]\pstgke[i]{n}{j} =
                     \dt{n}\resexpl[i](\pstgu{n}{j},\,
                                    \cpl[i](\pstgu[1]{n}{j},\,\dots,\,
                                            \pstgu[i-1]{n}{j},\,
                                            \pstp[i]{n-1},\,\dots,\,
                                            \pstp[m]{n-1},\,
                                            \mubold,\,\tstgi{n}{j}),\,
                                    \mubold,\,\tstgi{n}{j})}$
    \EndFor
   \EndFor
   \State Set
     $\displaystyle{\pstp{n} = \pstp{n-1} + \sum_{p=1}^s\hat{b}_p\pstgke{n}{p} +
                                            \sum_{p=1}^s b_p\pstgki{n}{p}}$
 \end{algorithmic}
\end{algorithm}
In this case, the Jacobian of the coupling predictor is block strictly lower
triangular
\begin{equation*}
 \pder{\cplprd}{\stvc} =
 \begin{bmatrix} 0 & & & \\
                 \displaystyle{\pder{\cpl[2]}{\stvc[1]}} & 0 & & \\
                 \vdots & \ddots & \ddots & \\
                 \displaystyle{\pder{\cpl[m]}{\stvc[1]}} & \cdots &
                 \displaystyle{\pder{\cpl[m]}{\stvc[m-1]}} & 0
 \end{bmatrix},
\end{equation*}
which implies the Jacobian of the monolithic implicit system is
block lower triangular
\begin{equation}
 D_{\stvc[j]} \resimpl[i] =
 \begin{cases}
  \displaystyle{\pder{\res[i]}{\stvc[i]}} & i = j \\
  \displaystyle{\pder{\res[i]}{\cpl[i]}
                               \pder{\cpl[i]}{\stvc[j]}} & i > j \\
  \zerobold & i < j.
 \end{cases}
 \label{EQ: JAC WEAK GS}
\end{equation}
This block lower triangular nature of the monolithic implicit system implies
that the individual systems can be solved sequentially beginning with system $1$
and yields a partitioned scheme. 

The implicit Jacobian of the monolithic implicit system of the weak
Gauss-Seidel predictor (\ref{EQ: JAC WEAK GS}) involves the entire lower
triangular portion of the coupling predictor; however, it is
{\it not required} for the implementation. From inspection of
\eqnref{EQ: IMEX STEP IMP EXP STAGE WITH PRED 3}, the implicit phase at stage $j$ for the $i$th
physical system requires the solution of a nonlinear system of equations in
the variable $\pstgu[i]{n}{j}$, with
$\pstgu[1]{n}{j},\,\dots,\,\pstgu[i-1]{n}{j}$ available from the implicit
solve corresponding to previous physical systems at the current stage.
Therefore, only the \emph{diagonal terms}
$\displaystyle{\frac{D \resimpl[i]}{D \stvc[i]} = \pder{\res[i]}{\stvc[i]}}$
of the monolithic implicit Jacobian are required, which shows that the Jacobians of
the coupling terms are not required for the weak Gauss-Seidel predictor.
This predictor is guaranteed to preserve the design order of the IMEX-RK
discretization and possesses stability properties in practice \cite{huang2018high}.

\section{Fully discrete sensitivity and adjoint method}\label{SEC: OPTIMIZATION}
In this section, we derive the expression for the total derivative of the quantity of interest $\mathcal{J}$ in \eqnref{EQ: QoI} with respect to the parameters $\mubold$, which is the essence in  gradient-based optimization.
 Since the evaluation of gradients is often the most costly step in the PDE-constraint optimization cycle, using efficient methods that accurately calculate the gradients are extremely important. There are generally two approaches to provide such information: the direct sensitivity approach and the adjoint approach \cite{giles2000introduction}. When the number of parameters  is smaller than the number of quantities of interest, the adjoint approach is much cheaper.

\subsection{Solver-consistent discretization of quantities of interest}
To maintain high-order accuracy for the optimization, discretization of the quantity of interest \eqnref{EQ: QoI} will be done in a 
solver-consistent manner \cite{zahr2016adjoint}, i.e. the spatial and temporal discretization used for the governing equation will also be 
used for the quantities of interest. The integral form \eqnref{EQ: QoI} can be rewritten as 
\begin{equation}
\label{EQ:QoI-SD}
 \pder{\Jcal}{t} = j(\stvc(t),\,\mubold,\,t).
\end{equation}
Augmenting the semi-discrete governing equations \eqnref{EQ: GOVERN-SD}(\ref{EQ: RES SPLIT}) with this ODE \eqnref{EQ:QoI-SD} yields the system of ODEs
\begin{equation} 
 \begin{bmatrix} \mass &  \\
                 & \mathbb{I}\\
 \end{bmatrix}
 \begin{bmatrix} \stvcdot  \\
                 \dot{\Jcal}\\
 \end{bmatrix}
 =
 \begin{bmatrix}  \resexpl(\stvc,\,\cplprd,\,\mubold,\,t)  \\
                  0\\
 \end{bmatrix}+
 \begin{bmatrix} \resimpl(\stvc,\,\cplprd,\,\mubold,\,t) \\
                 j(\stvc,\,\mubold,\,t)\\
 \end{bmatrix}.
\end{equation}

Applying the implicit-explicit temporal discretization introduced in Section~\ref{SEC: PARTITIONED SOLVER}  yields the fully discrete governing
equations and corresponding solver-consistent discretization of the quantity of interest \eqnref{EQ:QoI-SD}
\begin{equation*}
 \begin{aligned}
  \Jcal^0 &= 0, \\
  \Jcal^n &= \Jcal^{n-1}+
  \dt{n}\sum_{p=1}^s b_p  j(\pstgu{n}{p},\,\mubold,\,t_{n-1}+c_p\dt{n}).
 \end{aligned}
\end{equation*}

Finally, the objective functional in \eqnref{EQ: QoI} is evaluated at time
$t = T$ to yield the solver-consistent approximation 
\begin{equation}
\label{EQ: QoI-D}
 J(\pstgu{1}{1},\,\dots,\,\pstgu{n}{p},\,\mubold) = \Jcal^{N_t} =
 \sum_{n=1}^{N_t}\dt{n}
 \sum_{p=1}^s b_p j(\pstgu{n}{p},\,\mubold,\,t_{n-1}+c_p\dt{n}).
\end{equation}

\subsection{Direct sensitivity method}
Differentiation of the discretized weakly coupled Gauss-Seidel predictor based partitioned scheme expressions in Alg.~\ref{ALG: IMEX WEAK GS} with respect to $\mubold$ gives rise to the fully discrete sensitivity equations. For the $j$th stage  of the $n$th timestep, the sensitivity equations of the $i$the subsystem  write

\begin{subequations}
 \begin{align}
 \dmu[{\pstgu[i]{n}{j}}] &= \dmu[\ubm_{n-1}^i]  +
  \sum_{p=1}^{j-1}\hat{a}_{jp}\dmu[{\pstgke[i]{n}{p}}] +
  \sum_{p=1}^j a_{jp}\dmu[{\pstgki[i]{n}{p}}], \label{EQ: SENSITIVIT_1}\\
  \mass^i\dmu[{\pstgki[i]{n}{j}}] &=
  \dt{n}\left(\dmu[\gbm_{n,j}^{i}]+
              \pder{\gbm^{i}_{n,j}}{\pstgu[i]{n}{j}}\dmu[{\pstgu[i]{n}{j}}]+
              \pder{\gbm^i_{n,j}}{\cplprdstp[i]{n}}\dmu[{\cplprdstp[i]{n,j}}]\right), \label{EQ: SENSITIVIT_2}\\
  \mass^i\dmu[{\pstgke[i]{n}{j}}] &=
  \dt{n}\left(\dmu[\fbm_{n,j}^{i}]+
              \sum_{k=1}^{m}\pder{\fbm^{i}_{n,j}}{\pstgu[k]{n}{j}}\dmu[{\pstgu[k]{n}{j}}]+
              \pder{\fbm^i_{n,j}}{\cplprdstp[i]{n}}\dmu[{\cplprdstp[i]{n,j}}]\right), \label{EQ: SENSITIVIT_3}\\
  \dmu[\pstp{n}] &= \dmu[\pstp{n-1}] +
  \sum_{p=1}^s\hat{b}_p\dmu[\pstgke{n}{p}] +
  \sum_{p=1}^s b_p\dmu[\pstgki{n}{p}],\label{EQ: SENSITIVIT_4}
 \end{align}
\end{subequations}
here the $\cplprdstp[i]{n,j}$ is the weakly coupled Gaussian Seidel predictor in \eqnref{EQ: CPLPRD WEAK GS FD},  and its derivative with respect to $\mubold$ is 
\begin{equation}
  \dmu[{\cplprdstp[i]{n,j}}] = \dmu[{\cplstp[i]{n,j}}] + \sum_{p=1}^{i-1} \pder{\cplstp[i]{n,j}}{\pstgu[p]{n}{j}}\dmu[{\pstgu[p]{n}{j}}] + \sum_{p=i}^{m} \pder{\cplstp[i]{n,j}}{\ubm_{n-1}^p}\dmu[\ubm_{n-1}^p].
\end{equation}

By solving the sensitivities of the stage variables $\dmu[\pstgu{n}{p}]$ from \eqnref{EQ: SENSITIVIT_1}-\eqnref{EQ: SENSITIVIT_4}, the derivative of the quantity of interest, \eqnref{EQ: QoI-D}, of the multiphysics problem \eqnref{EQ: GOVERN} is written as
\begin{equation}\label{EQ: SENS DCOUPLE_DMU}
 \oder{J}{\mubold} =
 \sum_{n=1}^{N_t}\dt{n}
 \sum_{p=1}^s \hat{b}_p \pder{j(\pstgu{n}{p},\,t_{n-1}+c_p\dt{n})}{\pstgu{n}{p}} \dmu[\pstgu{n}{p}].
\end{equation}
Thanks to the partitioned nature of the multiphysics solver, the sensitivities of the stage variables $\dmu[\pstgu{n}{p}]$ can be solved substep-by-substep and subsystem-by-subsystem, the detailed algorithm is presented in Algorithm~\ref{ALG: IMEX WEAK GS SENS}.
\begin{algorithm}
 \caption{Direct sensitivity approach}
 \label{ALG: IMEX WEAK GS SENS}
 \begin{algorithmic}[1]
   \For{stages $j = 1,\,\dots,\,s$}
   \State Read stage solution $\pstp[i]{n-1}$, $\pstgke[i]{n}{p}$, $\pstgki[i]{n}{p}$ for $i = 1,\,\dots,\,m$ from disk.
      
    \For{physical systems $i = 1,\,\dots,\,m$}
     \State Construct $\dmu[{\cplprdstp[i]{n,j}}]$ based on \eqref{EQ: SENS DCOUPLE_DMU}
     \State Implicit solve \eqnref{EQ: SENSITIVIT_2} for $\dmu[{\pstgki[i]{n}{j}}]$:
     \State $\displaystyle \left(\mass^i - a_{jj}\pder{\gbm^{i}_{n,j}}{\pstgu[i]{n}{j}}\right)\dmu[{\pstgki[i]{n}{j}}] =
                                       \dt{n}\left(\dmu[\gbm_{n,j}^{i}]+
                                       \pder{\gbm^{i}_{n,j}}{\pstgu[i]{n}{j}}\left(\dmu[\ubm_{n-1}^i]  +
  \sum_{p=1}^{j-1}\hat{a}_{jp}\dmu[{\pstgke[i]{n}{p}}] +
  \sum_{p=1}^{j-1} a_{jp}\dmu[{\pstgki[i]{n}{p}}]\right)+
                                        \pder{\gbm^i_{n,j}}{\cplprdstp[i]{n}}\dmu[{\cplprdstp[i]{n,j}}]\right)$ 
    \State Construct \dmu[{\pstgu[i]{n}{j}}] based on \eqref{EQ: SENSITIVIT_1}
    \EndFor
    \For{physical systems $i = 1,\,\dots,\,m$}
     \State Explicit solve \eqnref{EQ: SENSITIVIT_3} for $\pstgke[i]{n}{j}$:  $\displaystyle{\dmu[{\pstgu[i]{n}{j}}] = \dmu[\ubm_{n-1}^i]  +
                                                                                 \sum_{p=1}^{j-1}\hat{a}_{jp}\dmu[{\pstgke[i]{n}{p}}] +
                                                                                 \sum_{p=1}^j a_{jp}\dmu[{\pstgki[i]{n}{p}}]}$
    \EndFor
   \EndFor
   \State Construct  $\dmu[\pstp{n}]$ based on \eqnref{EQ: SENSITIVIT_4}
 \end{algorithmic}
\end{algorithm} 
\subsection{Adjoint method}
The adjoint method provides an efficient alternative to the direct sensitivity method for evaluating the
 total derivative of the quantity of interest, especially when the number of parameters is large. 
 Before proceeding to the derivation of the adjoint equations, the following definitions are introduced for
the fully discrete Implicit-Explicit Runge-Kutta stage equations and state updates (See Alg.~\ref{ALG: IMEX WEAK GS})
\begin{equation}
\label{EQ: ADJ_PARA}
 \begin{aligned}
  \rstp[i]{0}(\pstp{0},\,\mubold) &= \pstp[i]{0} - \bar\ubm^{i}(\mubold), 
  \\
  \qstg[i]{n}{j}(\pstgu[i]{n}{j},\,\pstp[i]{n-1},\,\pstgke[i]{n}{1},\,\dots,\,\pstgke[i]{n}{j-1},\,
              \pstgki[i]{n}{1},\,\dots,\,\pstgki[i]{n}{j}) &=
             \pstgu[i]{n}{j} - \pstp[i]{n-1} -
             \sum_{p=1}^{j-1} \hat{a}_{jp}\pstgke[i]{n}{p} -
             \sum_{p=1}^j a_{jp}\pstgki[i]{n}{p}, 
  \\
  \rstgi[i]{n}{j}(\pstgu[i]{n}{j},\,\pstgki[i]{n}{j},\,\mubold,\,\cplprdstp[i]{n,j}) &=
  \mass\pstgki[i]{n}{j}-\dt{n}\gbm(\pstgu[i]{n}{j},\,\cplprdstp[i]{n,j},\,\mubold),
  \\
  \rstge[i]{n}{j}(\pstgu[i]{n}{j},\,\pstgke[i]{n}{j},\,\mubold,\,\cplprdstp[i]{n,j}) &=
  \mass\pstgke[i]{n}{j}-\dt{n}\fbm(\pstgu{n}{j},\,\cplprdstp[i]{n,j},\,\mubold), 
  \\
  \rstp[i]{n}(\pstp[i]{n},\,\pstp[i]{n-1},\,\pstgki[i]{n}{1},\,\dots,\,\pstgki[i]{n}{s},\,
  \pstgke[i]{n}{1},\,\dots,\,\pstgke[i]{n}{s}) &=
  \pstp[i]{n} - \pstp[i]{n-1} - \sum_{j=1}^s b_j \pstgki[i]{n}{j} -
  \sum_{j=1}^s \hat{b}_j \pstgke[i]{n}{j},
  \\
  \prstp[i]{n,j}(\pstgu[1]{n}{j},\,\dots,\,\pstgu[i-1]{n}{j},\,\pstp[i]{n-1} ,\,\dots,\,\pstp[m]{n-1},\,\cplprdstp[i]{n,j},\,\mubold) &=
              \cplprdstp[i]{n,j} - \cpl[i](\pstgu[1]{n}{j},\,\dots,\,
                                            \pstgu[i-1]{n}{j},\,
                                            \pstp[i]{n-1},\,\dots,\,
                                            \pstp[m]{n-1},\,
                                            \tstgi{n}{j},\,\mubold),
 \end{aligned}
\end{equation}
for $n = 1,\,\dots,\,N_t$, $i = 1,\,\dots,\,m$ and $j = 1,\,\dots,\,s$. Here $\bar\ubm^{i}(\mubold)$ is the initial condition, and in this work we use a steady-state solution to start the unsteady simulation.

Since the solution of the fully discretized PDE satisfies the above equations,
the QoI can be re-written as
\begin{equation*}
  J = J - \sum_{n=0}^{N_t}\sum_{i=1}^{m} \dstp{n}^i\rstp[i]{n}
        - \sum_{n=1}^{N_t}\sum_{j=1}^{s} \sum_{i=1}^{m} {\dstgi{n}{j}^i}\rstgi[i]{n}{j}
        - \sum_{n=1}^{N_t}\sum_{j=1}^{s} \sum_{i=1}^{m} {\dstge{n}{j}^i}\rstge[i]{n}{j}
        - \sum_{n=1}^{N_t}\sum_{j=1}^{s} \sum_{i=1}^{m} {\taustg{n}{j}^i}\qstg[i]{n}{j}
        - \sum_{n=1}^{N_t}\sum_{j=1}^{s} \sum_{i=1}^{m} {\sigstp{n,j}^i}\prstp[i]{n,j}
\end{equation*}
where $\dstp{n}^i$, ${\dstgi{n}{j}^i}$,  ${\dstge{n}{j}^i}$, ${\taustg{n}{j}^i}$, and ${\sigstp{n,j}^i}$ are test variables (also known as adjoint state variables or Lagrange multipliers) that
respectively enforce the state ODE system, coupling predictor, and initial conditions in \eqnref{EQ: ADJ_PARA}. Total differentiation of the modified QoI (or Lagrangian) leads to
\begin{footnotesize}
\begin{equation}
\label{EQ: ADJ_OPT_2}
 \begin{aligned}
  \oder{J}{\mubold} = \dmu[J] &+ \sum_{i=1}^m\dstp{0}^i\dmu[\bar\ubm^{i}] -
  \sum_{n=1}^{N_t}\sum_{i=1}^m \sum_{j=1}^s \dstgi{n}{j}^i\dmu[{\rstgi[i]{n}{j}}] -
  \sum_{n=1}^{N_t}\sum_{i=1}^m \sum_{j=1}^s \dstge{n}{j}^i\dmu[{\rstge[i]{n}{j}}] -
  \sum_{n=1}^{N_t}\sum_{i=1}^m \sum_{j=1}^s \sigstp{n,j}^i\dmu[{\prstp[i]{n,j}}] +
     \sum_{i=1}^m\left[-\dstp{N_t}^i\pder{\rstp[i]{N_t}}{\pstp[i]{N_t}}\right]\dmu[{\pstp[i]{N_t}}] 
  \\
  &+ \sum_{n=1}^{N_t}\sum_{i=1}^m \left[-\dstp{n-1}^i\pder{\rstp[i]{n-1}}{\pstp[i]{n-1}} -
           \dstp{n}^i\pder{\rstp[i]{n}}{\pstp[i]{n-1}} -
           \sum_{j=1}^s \taustg{n}{j}^i\pder{\qstg[i]{n}{j}}{\pstp[i]{n-1}} -
           \sum_{j=1}^s\sum_{p=1}^i\sigstp{n,j}^p\pder{\prstp[p]{n,j}}{\pstp[i]{n-1}}\right]\dmu[{\pstp[i]{n-1}}] 
  \\
  &+ \sum_{n=1}^{N_t}\sum_{i=1}^m\sum_{j=1}^s \left[-\dstp{n}^i\pder{\rstp[i]{n}}{\pstgki[i]{n}{j}}
           -\dstgi{n}{j}^i\pder{\rstgi[i]{n}{j}}{\pstgki[i]{n}{j}}
           -\sum_{p=j}^s\taustg{n}{p}^i\pder{\qstg[i]{n}{p}}{\pstgki[i]{n}{j}}
     \right]\dmu[{\pstgki[i]{n}{j}}] 
  \\
  &+ \sum_{n=1}^{N_t}\sum_{i=1}^m \sum_{j=1}^s\left[-\dstp{n}^i\pder{\rstp[i]{n}}{\pstgke[i]{n}{j}}
           -\dstge{n}{j}^i\pder{\rstge[i]{n}{j}}{\pstgke[i]{n}{j}}
           -\sum_{p=j+1}^s\taustg{n}{p}^i\pder{\qstg[i]{n}{p}}{\pstgke[i]{n}{j}}
     \right]\dmu[{\pstgke[i]{n}{j}}] 
  \\
  &+ \sum_{n=1}^{N_t}\sum_{i=1}^m \sum_{j=1}^s \left[\pder{J}{\pstgu[i]{n}{j}}
           -\dstgi{n}{j}^i\pder{\rstgi[i]{n}{j}}{\pstgu[i]{n}{j}}
           -\sum_{k=1}^{m}\dstge{n}{j}^{k}\pder{\rstge[k]{n}{j}}{\pstgu[i]{n}{j}}
           -\taustg{n}{j}^i\pder{\qstg[i]{n}{j}}{\pstgu[i]{n}{j}}
           -\sum_{p = i+1}^{m}\sigstp{n,j}^p\pder{\prstp[p]{n,j}}{\pstgu[i]{n}{j}}
     \right]\dmu[{\pstgu[i]{n}{j}}]
  \\
  &+ \sum_{n=1}^{N_t}\sum_{i=1}^m \sum_{j=1}^s \left[
           -\dstgi{n}{j}^i\pder{\rstgi[i]{n}{j}}{\cplprdstp[i]{n,j}}
           -\dstge{n}{j}^i\pder{\rstge[i]{n}{j}}{\cplprdstp[i]{n,j}}
           -\sigstp{n,j}^i\pder{\prstp[i]{n,j}}{\cplprdstp[i]{n,j}}
  \right]\dmu[{\cplprdstp[i]{n,j}}],
 \end{aligned}
\end{equation}
\end{footnotesize}
here, we re-arrange these terms, such that the state variable sensitivities are isolated. The adjoint state variables 
$\dstp{n}^i$, ${\dstgi{n}{j}^i}$,  ${\dstge{n}{j}^i}$, ${\taustg{n}{j}^i}$, and ${\sigstp{n,j}^i}$,
 which have remained arbitrary to this point, are chosen such that the
bracketed terms in \eqnref{EQ: ADJ_OPT_2} vanish. The adjoint equations are 
\begin{subequations}
 \begin{align}
  \dstp{N_t}^i &= \zerobold, \label{EQ: ADJ 1}
  \\
  \dstp{n-1}^i &= \dstp{n}^i + \sum_{j=1}^s \taustg{n}{j}^i +
           \sum_{j=1}^s\sum_{p=1}^i\pder{\cplstp[p]{n,j}}{\pstp[i]{n-1}}^{T}\sigstp{n,j}^p,\label{EQ: ADJ 2}
  \\
  {\mass^{i}}^{T}\dstgi{n}{j}^i &= b_j\dstp{n}^i + \sum_{p=j}^s a_{pj}\taustg{n}{p}^i, \label{EQ: ADJ 3}
  \\
  {\mass^{i}}^{T}\dstge{n}{j}^i &= \hat{b}_j\dstp{n}^i + \sum_{p=j+1}^s \hat{a}_{pj}\taustg{n}{p}^i, \label{EQ: ADJ 4}
  \\
  \taustg{n}{j}^{i} &= \pder{J}{\pstgu[i]{n}{j}} +
                    \dt{n}{\pder{\gbm_{n,j}^i}{\pstgu[i]{n}{j}}}^T\dstgi{n}{j}^i +
                    \dt{n}\sum_{k=1}^{m}{\pder{\fbm_{n,j}^k}{\pstgu[i]{n}{j}}}^T\dstge{n}{j}^k +
                    \sum_{p = i+1}^{m}\pder{\cplstp[p]{n,j}}{\pstgu[i]{n}{j}}^{T}\sigstp{n,j}^p,\label{EQ: ADJ 5}
  \\
  \sigstp{n,j}^{i} &=
           \dt{n}\pder{\gbm_{n,j}^i}{\cplprdstp[i]{n,j}}^T\dstgi{n}{j}^{i} +
           \dt{n}\pder{\fbm_{n,j}^i}{\cplprdstp[i]{n,j}}^T\dstge{n}{j}^{i},\label{EQ: ADJ 6}
 \end{align}
\end{subequations}
for $n = 1,\,\dots,\,N_t$, $i = 1,\,\dots,\,m$ and $j = 1,\,\dots,\,s$. These are the fully discrete adjoint equations corresponding to the multiphysics problem
in \eqnref{EQ: GOVERN}, discrete quantity of interest $\mathcal{J}$, and parameter $\mubold$. Solving the
adjoint variables reversely from \eqnref{EQ: ADJ 1}-\eqnref{EQ: ADJ 6}, the expression for the gradient in \eqnref{EQ: ADJ_OPT_2} reduces to
\begin{equation*}
 \begin{aligned}
  \oder{J}{\mubold} = \dmu[J] &+ \sum_{i=1}^m\dstp{0}^i\dmu[\bar\ubm^{i}] +
  \sum_{n=1}^{N_t}\sum_{i=1}^m \sum_{j=1}^s \dt{n}\pder{\gbm_{n,j}^i}{\mu}^T\dstgi{n}{j}^{i} +
  \sum_{n=1}^{N_t}\sum_{i=1}^m \sum_{j=1}^s \dt{n}\pder{\fbm_{n,j}^i}{\mu}^T\dstge{n}{j}^{i} +
  \sum_{n=1}^{N_t}\sum_{i=1}^m \sum_{j=1}^s \dmu[{\cplstp[i]{n,j}}]^T \sigstp{n,j}^i
   \end{aligned}
\end{equation*}

Due to the partitioned nature of the multiphysics solver, the adjoint variables can be solved substep-by-substep and subsystem-by-subsystem, the detailed algorithm is presented in Algorithm~\ref{ALG: IMEX WEAK GS ADJOINT}.
\begin{algorithm}
 \caption{Adjoint approach}
 \label{ALG: IMEX WEAK GS ADJOINT}
 \begin{algorithmic}[1]
   \For{stages $j = s,\,\dots,\,1$}
   \State Read stage solution $\pstp[i]{n-1}$, $\pstgke[i]{n}{p}$, $\pstgki[i]{n}{p}$ for $i = 1,\,\dots,\,m$ from disk.
      
    \For{physical systems $i = m,\,\dots,\,1$}
     \State Explicit solve \eqnref{EQ: ADJ 4} for $\dstge{n}{j}^i$:
            $\displaystyle {\mass^{i}}^{T}\dstge{n}{j}^i = \hat{b}_j\dstp{n}^i + \sum_{p=j+1}^s \hat{a}_{pj}\taustg{n}{p}^i $ 
    \EndFor
    \For{physical systems $i = m,\,\dots,\,1$}
     \State Set $\displaystyle \widetilde{\taustg{n}{j}^{i}} = \pder{J}{\pstgu[i]{n}{j}} +
                    \dt{n}\sum_{k=1}^{m}{\pder{\fbm_{n,j}^k}{\pstgu[i]{n}{j}}}^T\dstge{n}{j}^k +
                    \sum_{p = i+1}^{m}\pder{\cplstp[p]{n,j}}{\pstgu[i]{n}{j}}^{T}\sigstp{n,j}^p$
     \State Implicit solve \eqnref{EQ: ADJ 3} for $\dstgi{n}{j}^i$:  $\displaystyle \left({\mass^{i}}^{T} -\dt{n}{\pder{\gbm_{n,j}^i}{\pstgu[i]{n}{j}}}^T\right) \dstgi{n}{j}^i = b_j\dstp{n}^i + \sum_{p=j+1}^s a_{pj}\taustg{n}{p}^i + a_{jj}\widetilde{\taustg{n}{j}^{i}}$
     \State Construct $\taustg{n}{j}^{i}$ and $\sigstp{n,j}^{i}$ based on \eqref{EQ: ADJ 5} and \eqref{EQ: ADJ 6}
    \EndFor
   \EndFor
   \State Construct  $\dstp{n-1}^i$ based on \eqref{EQ: ADJ 2}
 \end{algorithmic}
\end{algorithm} 
\section{Applications} \label{SEC: APPLICATIONS}
In this section, we demonstrate the proposed high-order optimization procedure on two multiphysics problems: 
a 1D fluid-structure-mesh three-field coupling piston problem and a 2D fluid-structure two-field coupling foil energy harvesting problem.
\subsection{Governing equations and semi-discretization}
\subsubsection{Compressible fluid flow}
\label{SUBSUBSEC: ALE Fluid}
\begin{figure}
  \centering
  \includegraphics[width=2.5in]{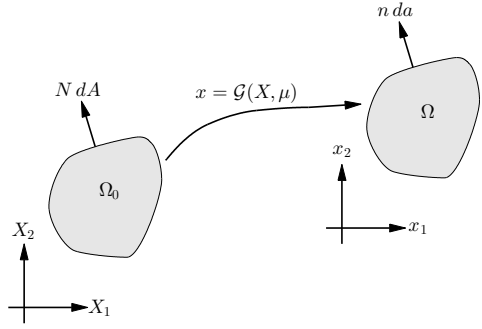}
  \caption{Mapping between reference and physical domains.}
  \label{FIG: DOM MAP}
\end{figure}

The governing equations for compressible fluid flow, defined on a
deformable fluid domain $\Omega(\mubold,\,t)$, can be written as a viscous conservation
law
\begin{equation}
\label{EQ: FSI GOVERN}
\frac{\partial U}{\partial t} + \nabla \cdot \Fcal^{inv}(U) + \nabla \cdot \Fcal^{vis}(U, \nabla U)= 0 \quad \text{in} \quad \Omega(\mubold, t),
\end{equation}
where $U$ is the conservative state variable vector and the physical
flux consists of an inviscid part $\Fcal^{inv}(U)$ and a viscous part
$\Fcal^{vis}(U,\,\nabla U)$. The conservation law in (\ref{EQ: FSI GOVERN})
is transformed to a fixed reference domain $\Omega_0$ by defining a
time-dependent diffeomorphism $\Gcal$ between the reference domain and
the physical domain; see Figure~\ref{FIG: DOM MAP}. At each time $t$, a point
$X$ in the reference domain $\Omega_0$ is mapped to $x(X,\mubold, t) = \Gcal(X,\mubold,t)$
in the physical domain $\Omega(\mubold,t)$.

The deformation gradient $G$, velocity $v_G$, and Jacobian $g$ of
the mapping are defined as
\begin{equation}
 G = \nabla_X\Gcal\,,\quad
 v_G = \pder{\Gcal}{t}\,,\quad
 g = \det G.
\end{equation}
Following the procedure in \cite{persson2009discontinuous, zahr2016adjoint},
the governing equation
(\ref{EQ: FSI GOVERN}) can be written in the reference domain as
\begin{equation}\label{EQ: FSI GOVERN TRANSF}
\frac{\partial U_X}{\partial t} + \nabla_X \cdot \Fcal_X^{inv}(U_X) +
\nabla_X\cdot \Fcal_X^{vis}(U_X, \nabla_X U_X)= 0
\quad \text{in} \quad \Omega_0, 
\end{equation}
where $\nabla_X$ defines the spatial derivative with respect to the reference
domain, conserved quantities and its derivatives in the reference domain are
written as
\begin{equation}\label {EQ: FSI GOVERN TRANSF2}
 U_X = gU\,,\qquad
 \nabla_X U_X = g\nabla U_X \cdot G + g^{-1}U_X \frac{\partial g}{\partial X}.
\end{equation}
The inviscid and viscous fluxes are transformed to the reference domain as
\begin{equation}
\begin{aligned}
 \Fcal_X^{inv}(U_X) &= g\Fcal^{inv}(g^{-1} U_X)G^{-T} -
                       U_X \mathbin{\mathop{\otimes}} G^{-1}v_G, \\
 \Fcal_X^{vis}(U_X) &= g\Fcal^{vis}\left(g^{-1} U_X,
                       g^{-1}\left[\nabla_X U_X -
                       g^{-1} U_X \pder{g}{X}\right]G^{-1}\right)G^{-T}.
 \end{aligned}
\end{equation}

The governing equations in (\ref{EQ: FSI GOVERN TRANSF}) reduce to the following
system of ODEs after an appropriate spatial discretization, such as a
discontinuous Galerkin or finite volume method, is applied
\begin{equation}
\label{EQ: FSI3 FLUID COUPLING}
 \mass[f]\stvcdot[f] = \res[f](\stvc[f],\,\cpl[f],\,\mubold,\,t),
\end{equation}
where $\mass[f]$ is the fixed mass matrix, $\stvc[f]$ is the semi-discrete
fluid state vector, i.e., the discretization of $U_X$ on $\Omega_0$,
$\res[f](\stvc[f],\,\cpl[f],\,\mubold,\,t)$ is the spatial discretization of the transformed
inviscid and viscous fluxes on $\Omega_0$, and $\cpl[f]$ is the coupling term
that might contain information about the domain mapping $\Gcal(X,\,\mubold,\,t)$. In
particular, the coupling term contains the position and velocities of the
nodal coordinates of the computational mesh.
The domain mapping is defined using an element-wise nodal (Lagrangian)
polynomial basis on the mesh with coefficients from the nodal positions and
velocities.

\subsubsection{Simple structure model}
\label{SUBSUBSEC: APP FSI STRUCT}
In general, the governing equations for the structure will be given by
a system of partial differential equation such as the continuum equations in
total Lagrangian form with an arbitrary constitutive law. However, in this
work, we only consider simple structures like mass-spring-damper systems
that can directly be written as a second-order system of ODEs
\begin{equation} \label{EQ: SIMP STRUCT}
 m_s\ddot{u}_s + c_s\dot{u}_s + k_s u_s = f_{ext}(t),
\end{equation}
where $m_s$ is the mass of the (rigid) object, $c_s$ is the damper resistance
constant, $k_s$ is the spring stiffness, and $f_{ext}(t)$ is a time-dependent
external load, which will be given by integrating the pointwise force the
fluid exerts on the object. 

The equations in
(\ref{EQ: SIMP STRUCT}) are re-written in a first-order form, to conform to the notation in this document, as
\begin{equation}
 \mass[s]\stvcdot[s] = \res[s](\stvc[s],\,\cpl[s],\,\mubold,\,t).
\end{equation}
In the case of the simple structure in (\ref{EQ: SIMP STRUCT}), the mass
matrix, state vector, residual, and coupling term are
\begin{equation}
 \mass[s] = \begin{bmatrix} m_s & \\ & 1 \end{bmatrix}, \qquad
 \stvc[s] = \begin{bmatrix} \dot{u}_s \\ u_s \end{bmatrix}, \qquad
 \cpl[s] = f_{ext}, \qquad
 \res[s](\stvc[s],\,\cpl[s]) =
              \begin{bmatrix} f_{ext}-c_s\dot{u}_s-k_s u_s \\ u_s \end{bmatrix}.
\label{EQ: FSI3 STRUCT COUPLING}
\end{equation}

\subsubsection{Deformation of the fluid domain}
The mesh deformation is generally described by a pseudo-structure driven solely
by Dirichlet boundary conditions provided by the displacement of the structure
at the fluid-structure interface \cite{farhat1995mixed, farhat1998torsional} or a
parametrized mapping such as radial basis functions \cite{van2007higher, rendall2008unified, froehle2014high} or
blending maps \cite{persson2009discontinuous}. Due to different treatments of the mesh deformation, the fluid-structure interaction problem can be formulated as 
three-field coupling or two-filed coupling problems.

For the first treatment, the governing equations are given by the
continuum mechanics equations in total Lagrangian form with an arbitrary
constitutive law
\begin{equation} 
 \begin{aligned}
  \pder{\bar{p}}{t} - \nabla \cdot P(G) &= 0
    \qquad\qquad &&\text{in}~\Omega_0 \\
  x &= x_b
    \qquad\qquad &&\text{on}~\partial \Omega_0^D \\
  \dot{x} &= \dot{x}_b
    \qquad\qquad &&\text{on}~\partial \Omega_0^D,
 \end{aligned}
\label{EQ: FSI3 MESH}
\end{equation}
where $\bar{p}(X,\,t) = \rho_m\dot{x}$ is the linear momentum, $\rho_m$
is the density, and $P$ is the first Piola-Kirchhoff stress of the
pseudo-structure. The deformation gradient $G$ is the mapping that
defines the deformation of the reference fluid domain $\Omega_0$ to physical
fluid domain $\Omega(t)$. The position and velocity of the fluid domain are
prescribed along $\partial \Omega_0^D$, the union of the fluid-structure
interface and the fluid domain boundary.
The governing equations in (\ref{EQ: FSI3 MESH}) reduce to the following
system of ODEs after an appropriate spatial discretization, such as the
finite element method, is applied and recast in first-order form
\begin{equation}
 \Mbm^x\dot\ubm^x = \rbm^x(\ubm^x,\,\cbm^x,\,\mubold,\,t)
\end{equation}
where $\mass[x]$ is the fixed mass matrix, $\stvc[x](t)$ is the semi-discrete
state vector consisting of the displacements and velocities of the
mesh nodes, $\res[x](\stvc[x],\,\cpl[x],\,\mubold,\,t)$ is the spatial discretization
of the continuum equations and boundary conditions on the reference domain
$\Omega_0$, and $\cpl[x]$ is the coupling term that contains information about
the motion of the fluid structure interface. This model of the mesh motion
leads to a three-field FSI formulation when coupled to the fluid and
structure equations.

For the second treatment, the domain mapping
$x = \Gcal(X,\,t)$ is given by an analytical function, parametrized by the
deformation and velocity of the fluid-structure interface, that can be
analytically differentiated to obtain the deformation gradient $G(X,\,t)$
and velocity $v_G(X,\,t)$. Since the fluid mesh motion is no longer included
in the system of time-dependent partial differential equations, this leads to
a two-field FSI formulation in terms of the fluid and structure states only.

\subsubsection{Two-field and three-field fluid-structure coupling}
In the three-field fluid-structure interaction setting
\begin{equation}
\label{EQ: FSI3 SEMI DISC}
\mass[s]\stvcdot[s] = \res[s](\stvc[s],\,\cpl[s],\,\mubold,\,t), \quad 
\mass[x]\stvcdot[x] = \res[x](\stvc[x],\,\cpl[x],\,\mubold,\,t) , \quad
\mass[f]\stvcdot[f] = \res[f](\stvc[f],\,\cpl[f],\,\mubold,\,t)
\end{equation}
introduced in \cite{farhat1995mixed}, the coupling terms have the following
dependencies
\begin{equation} 
\cpl[s] = \cpl[s](\stvc[s],\,\stvc[x],\,\stvc[f],\,\mubold,\,t), \quad
\cpl[x] = \cpl[x](\stvc[s],\,\mubold,\,t), \quad
\cpl[f] = \cpl[f](\stvc[s],\, \stvc[x],\,\mubold,\,t).
\label{EQ: FSI3 COUPLE}
\end{equation}
From \eqnref{EQ: FSI3 STRUCT COUPLING}, the structure coupling term is the external force applied
to the structure that comes from integrating the fluid stresses over the
fluid-structure interface. The mesh coupling term is the position and
velocity of the fluid-structure interface and therefore depends solely on
the state of the structure. From Eq.~(\ref{EQ: FSI GOVERN TRANSF})-(\ref{EQ: FSI GOVERN TRANSF2}), 
the fluid coupling term is
the position and velocity of the entire fluid mesh and therefore depends
on the state of the structure and the mesh.

In the two-field FSI setting
\begin{equation}
\label{EQ: FSI2 SEMI DISC}
\mass[s]\stvcdot[s] = \res[s](\stvc[s],\,\cpl[s],\,\mubold,\,t), \quad 
\mass[f]\stvcdot[f] = \res[f](\stvc[f],\,\cpl[f],\,\mubold,\,t)
\end{equation}
the mesh motion is given by an analytical function and the coupling terms have
the following dependencies
\begin{equation} \label{EQ: FSI2 COUPLE}
\cpl[s] = \cpl[s](\stvc[s],\,\stvc[f],\,\mubold,\,t), \quad
\cpl[f] = \cpl[f](\stvc[s],\,\mubold,\,t).
\end{equation}
In this case, the structure coupling term is determined from the fluid and
structure state since the external force depends on the traction integrated
over the fluid-structure interface. The fluid coupling term, i.e., the
position and velocity of the fluid mesh, is determined from the structure
state. Finally, the ordering of the subsystems implied in
(\ref{EQ: FSI3 SEMI DISC}) and (\ref{EQ: FSI2 SEMI DISC}) is used throughout
the remainder of this section, which plays an important role when
defining the Gauss-Seidel predictors.

\subsection{1D fluid-structure-mesh three-field coupling piston problem}
This proposed optimization procedure is first verified by the canonical FSI
model problem: a one-dimensional piston problem (\figref{FIG: 1D PISTON}).
\begin{figure}
\centering
    \begin{tikzpicture}[scale=1.3]
      \fill [pattern = north east lines] (-0.3,0) rectangle (0.3,1.0);
      \fill[cyan] (0,0) rectangle (4,1.0);
      \draw [fill=red!15!white,line width=0.5mm] (4,0) rectangle (4.5,1);
      \draw node at (1.5,0.5) {inviscid flow};
      
      \draw[decoration={aspect=0.3, segment length=1.5mm, amplitude=3mm,coil},decorate] (4.5,0.5) -- (5.5,0.5);
      \fill [pattern = north east lines] (5.5,0) rectangle (6.0,1.0);
      \draw[line width=0.5mm](0, 0)-- (0, 1);
      \draw[line width=0.5mm](5.5, 0)-- (5.5, 1);
      \draw[thick,->,line width=0.5mm,](0,-0.3) -- (4.5,-0.3) node[anchor=north west][scale=1.2] {x};
      \draw node at (4.25,0.5) {$m$};
      \draw node at (5.0,1) {$k$};
    \end{tikzpicture}
    \caption{One-dimensional piston system}
    \label{FIG: 1D PISTON}
\end{figure}
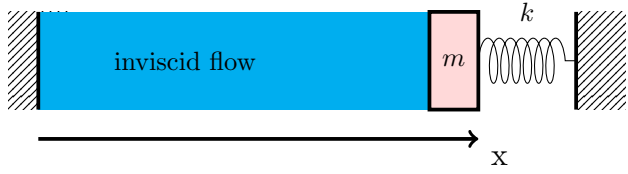
The inviscid fluid is governed by the one-dimensional Euler equations defined on 
$x \in \Omega(t) = [0,\,1.0-u_s]$, where $u_s$ is the displacement of the
piston. The fluid flow is the adiabatic gas with constant $\gamma = 1.4$. The fluid is
initially at rest $u = 0$ with a density $\rho = 1.0$ and pressure $p = 0.4$.
After transformation to the reference domain $\Omega_0 = [0,\,1]$ following the
procedure in Section~\ref{SUBSUBSEC: ALE Fluid}, the equations are semi-discretized by a standard
first-order finite volume method using Roe's flux \cite{roe1981approximate}
with $100$ elements.

The deformation of the fluid mesh is handled by considering the fluid domain
to be a pseudo-structure governed by the continuum equations in
\eqnref{EQ: FSI3 MESH}, restricted to the one-dimensional case with a
linear constitutive law and infinitesimal strains assumed
\begin{equation}
\rho_m\ddot{u}_x = E_m\frac{\partial^2 u_x}{\partial X^2} - c_m \dot{u}_x,
\end{equation}
where $u_x(X,\,t)$ is the mesh displacement vector defined over the reference
domain $X \in \Omega_0$ and the density, Young's modulus, and damping
coefficient are $\rho_m = 1.0$, $E_m = 1.0$, $c_m = 0.0$, respectively. The
governing equation for the mesh deformation is discretized in space using the
finite difference method.

Finally, the structure is modeled by a linear mass-spring system as
\eqnref{EQ: SIMP STRUCT} with piston mass $m_s = 1.0$, spring stiffness
$\mu_k = 1.0$, and no damper $c_s = 0$. The piston is initially displaced a
distance of $u_s = 0.0$. Once the piston is released, it immediately begins
to recede due to the combination of the spring being perturbed from its
equilibrium configuration and the flow pressure, which causes a $C^0$
rarefaction wave near the interface.

The objective function to minimize is set to be the integral of square of the piston displacement till $T=1.0$
\begin{equation}
 \mathcal{J} = \int_{0}^{T} u_s^2 dt.
\end{equation}
The only parameter is the stiffness of the piston $\mu_k$, for verification purpose, an additional constraint $0 \leq \mu_k \leq 10$ is imposed. We should expect that when the stiffness reaches its maximum, the objective function reaches its minimum.

Table~\ref{TABLE: 1D FSI DERIVATIVES} shows the objective function and its derivative evaluated in three different ways, using central finite differences with $\epsilon = 10^{-6}$, the direct sensitivity method, and the adjoint method.
The results of the direct sensitivity method and the adjoint method are within $10^{-6}$ of the finite difference results, which verifies the correctness of our current implementation. Moreover, the accuracy of the finite difference method is limited by the ``step-size dilemma,'' therefore the adjoint method and the direct sensitivity method are likely producing more accurate derivatives. 

\begin{table}[!htb]
\begin{tabular}{c|cccc}
Scheme & $\mathcal{J}$ &  FD  &  Direct &   Adjoint\\
\hline
IMEX1&   5.24027644581e-03   &- 6.40416043546e-04     & - 6.40416045418e-04      & -6.40416045418e-04 \\
IMEX2&   5.01357571586e-03   &- 5.75379291520e-04     &  -5.75379340362e-04      & -5.75379340362e-04 \\
IMEX3&   5.01291619482e-03    &-5.75053709945e-04     &  -5.75053861151e-04      & -5.75053861151e-04  \\
IMEX4&   5.01291415604e-03    &-5.75054676186e-04     &  -5.75054797593e-04      & -5.75054797593e-04  \\
\end{tabular}
\caption{1D piston problem: the objective function value and its gradients.}\label{TABLE: 1D FSI DERIVATIVES}
\end{table}

The convergence of the quantities of interest is reported in \figref{FIG: 1D PISTON}-left, the corresponding convergence of the parameter is reported in \figref{FIG: 1D PISTON}-right. All IMEX schemes use step size $\Delta t = 0.01$ and lead to convergence in 8 optimization steps. The parameter $\mu_k$ converges to its upper bound as expected. 
\begin{figure}
\begin{tikzpicture}
\begin{axis}[
    width=0.48\textwidth,
    height=0.48\textwidth,
    ylabel={$\mathcal{J}$},
    xlabel={Iteration}]
\addplot [red, solid, thick, mark=square*, mark size=1, mark options={solid}]  table[x index=0, y index=1] {data/piston_euler/euler_1d_3fields_optim_RK1.dat};\label{line:1d_fsi:rk1:obj}
\addplot [gray, loosely dotted, thick, mark=otimes*, mark size=1, mark options={solid}]  table[x index=0, y index=1] {data/piston_euler/euler_1d_3fields_optim_RK2.dat};\label{line:1d_fsi:rk2:obj}
\addplot [green!75!black, dashed, thick, mark=triangle*, mark size=1, mark options={solid}]  table[x index=0, y index=1] {data/piston_euler/euler_1d_3fields_optim_RK3.dat};\label{line:1d_fsi:rk3:obj}
\addplot [blue, dashdotted, thick, mark=diamond*, mark size=1, mark options={solid}]  table[x index=0, y index=1] {data/piston_euler/euler_1d_3fields_optim_RK4.dat};\label{line:1d_fsi:rk4:obj}
\end{axis}

\end{tikzpicture}\quad 
\begin{tikzpicture}
\begin{axis}[
    width=0.48\textwidth,
    height=0.48\textwidth,
    ylabel={$\mu_k$},
    xlabel={Iteration}]
\addplot [red, solid, thick, mark=square*, mark size=1, mark options={solid}]  table[x index=0, y index=2] {data/piston_euler/euler_1d_3fields_optim_RK1.dat};
\addplot [gray, loosely dotted, thick, mark=otimes*, mark size=1, mark options={solid}]  table[x index=0, y index=2] {data/piston_euler/euler_1d_3fields_optim_RK2.dat};
\addplot [green!75!black, dashed, thick, mark=triangle*, mark size=1, mark options={solid}]  table[x index=0, y index=2] {data/piston_euler/euler_1d_3fields_optim_RK3.dat};
\addplot [blue, dashdotted, thick, mark=diamond*, mark size=1, mark options={solid}]  table[x index=0, y index=2] {data/piston_euler/euler_1d_3fields_optim_RK4.dat};
\end{axis}

\end{tikzpicture}
\caption{Convergence of the IMEX1 (\ref{line:1d_fsi:rk1:obj}),  IMEX2 (\ref{line:1d_fsi:rk2:obj}), IMEX3 (\ref{line:1d_fsi:rk3:obj}), and IMEX4 (\ref{line:1d_fsi:rk4:obj}) schemes when applied to the three-field coupling piston problem: the value of the objective function at each iteration (left); the value of the parameter at each iteration (right).}
\label{FIG: 1D PISTON}
\end{figure}
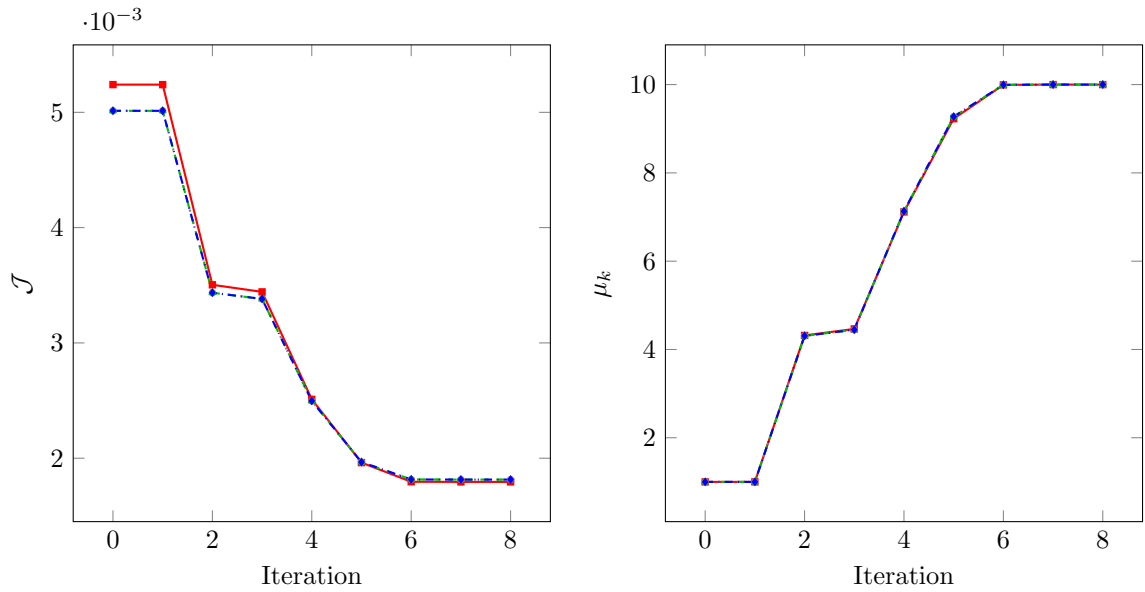

\subsection{2D fluid-structure two-field coupling foil energy harvesting process}

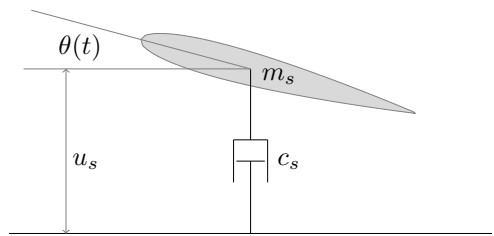
\begin{figure}
 \centering
 \begin{tikzpicture}
\begin{axis}[
axis equal,
axis lines=none,
width=8cm,
ymax=0.6,
xmax=0.85,
xmin=-0.85,
ymin=-0.6]
\addplot [white, forget plot]
coordinates {
( -0.85000000,  -0.60000000)
(  0.85000000,  -0.60000000)
(  0.85000000,   0.60000000)
( -0.85000000,   0.60000000)
( -0.85000000,  -0.60000000)};

\draw [smooth, ultra thin] plot coordinates {(axis cs:-1.000000, -0.580000) (axis cs:1.000000, -0.580000)
};

\draw [smooth, ultra thin, opacity=0.5, fill=black!30!white] plot coordinates {(axis cs:-0.385791, 0.103372) (axis cs:-0.371603, 0.117293) (axis cs:-0.360141, 0.121044) (axis cs:-0.349136, 0.123087) (axis cs:-0.338373, 0.124230) (axis cs:-0.327767, 0.124784) (axis cs:-0.317276, 0.124909) (axis cs:-0.306874, 0.124701) (axis cs:-0.296545, 0.124224) (axis cs:-0.286276, 0.123522) (axis cs:-0.276059, 0.122625) (axis cs:-0.265887, 0.121559) (axis cs:-0.255756, 0.120343) (axis cs:-0.245660, 0.118992) (axis cs:-0.235597, 0.117520) (axis cs:-0.225564, 0.115935) (axis cs:-0.215559, 0.114249) (axis cs:-0.205579, 0.112467) (axis cs:-0.195623, 0.110598) (axis cs:-0.185688, 0.108646) (axis cs:-0.175775, 0.106617) (axis cs:-0.165880, 0.104516) (axis cs:-0.156004, 0.102346) (axis cs:-0.146146, 0.100111) (axis cs:-0.136304, 0.097816) (axis cs:-0.126477, 0.095462) (axis cs:-0.116665, 0.093053) (axis cs:-0.106868, 0.090591) (axis cs:-0.097083, 0.088079) (axis cs:-0.087312, 0.085518) (axis cs:-0.077553, 0.082912) (axis cs:-0.067806, 0.080262) (axis cs:-0.058070, 0.077570) (axis cs:-0.048345, 0.074838) (axis cs:-0.038630, 0.072066) (axis cs:-0.028925, 0.069258) (axis cs:-0.019230, 0.066413) (axis cs:-0.009544, 0.063535) (axis cs:0.000133, 0.060623) (axis cs:0.009802, 0.057679) (axis cs:0.019462, 0.054704) (axis cs:0.029114, 0.051700) (axis cs:0.038759, 0.048667) (axis cs:0.048396, 0.045606) (axis cs:0.058026, 0.042519) (axis cs:0.067649, 0.039405) (axis cs:0.077266, 0.036267) (axis cs:0.086876, 0.033105) (axis cs:0.096480, 0.029920) (axis cs:0.106077, 0.026711) (axis cs:0.115669, 0.023481) (axis cs:0.125255, 0.020230) (axis cs:0.134836, 0.016958) (axis cs:0.144411, 0.013666) (axis cs:0.153981, 0.010355) (axis cs:0.163546, 0.007025) (axis cs:0.173106, 0.003676) (axis cs:0.182662, 0.000310) (axis cs:0.192212, -0.003074) (axis cs:0.201758, -0.006475) (axis cs:0.211300, -0.009893) (axis cs:0.220837, -0.013327) (axis cs:0.230370, -0.016776) (axis cs:0.239899, -0.020241) (axis cs:0.249424, -0.023722) (axis cs:0.258945, -0.027217) (axis cs:0.268462, -0.030727) (axis cs:0.277975, -0.034251) (axis cs:0.287484, -0.037789) (axis cs:0.296989, -0.041342) (axis cs:0.306491, -0.044908) (axis cs:0.315989, -0.048488) (axis cs:0.325484, -0.052082) (axis cs:0.334974, -0.055688) (axis cs:0.344462, -0.059309) (axis cs:0.353946, -0.062942) (axis cs:0.363426, -0.066589) (axis cs:0.372903, -0.070248) (axis cs:0.382376, -0.073921) (axis cs:0.391845, -0.077607) (axis cs:0.401312, -0.081306) (axis cs:0.410774, -0.085018) (axis cs:0.420234, -0.088743) (axis cs:0.429689, -0.092482) (axis cs:0.439141, -0.096234) (axis cs:0.448590, -0.099999) (axis cs:0.458034, -0.103778) (axis cs:0.467475, -0.107571) (axis cs:0.476913, -0.111377) (axis cs:0.486346, -0.115198) (axis cs:0.495776, -0.119033) (axis cs:0.505202, -0.122882) (axis cs:0.514624, -0.126747) (axis cs:0.524042, -0.130626) (axis cs:0.533456, -0.134520) (axis cs:0.542865, -0.138430) (axis cs:0.552271, -0.142355) (axis cs:0.561672, -0.146297) (axis cs:0.571069, -0.150255) (axis cs:0.579809, -0.156664) (axis cs:0.569688, -0.155410) (axis cs:0.559571, -0.154139) (axis cs:0.549458, -0.152852) (axis cs:0.539350, -0.151549) (axis cs:0.529246, -0.150230) (axis cs:0.519146, -0.148896) (axis cs:0.509051, -0.147546) (axis cs:0.498959, -0.146182) (axis cs:0.488871, -0.144802) (axis cs:0.478787, -0.143409) (axis cs:0.468707, -0.142001) (axis cs:0.458631, -0.140579) (axis cs:0.448558, -0.139143) (axis cs:0.438489, -0.137693) (axis cs:0.428424, -0.136230) (axis cs:0.418363, -0.134753) (axis cs:0.408304, -0.133263) (axis cs:0.398250, -0.131760) (axis cs:0.388199, -0.130243) (axis cs:0.378152, -0.128713) (axis cs:0.368108, -0.127170) (axis cs:0.358067, -0.125615) (axis cs:0.348030, -0.124046) (axis cs:0.337997, -0.122463) (axis cs:0.327967, -0.120868) (axis cs:0.317941, -0.119260) (axis cs:0.307918, -0.117638) (axis cs:0.297899, -0.116003) (axis cs:0.287883, -0.114354) (axis cs:0.277871, -0.112691) (axis cs:0.267863, -0.111015) (axis cs:0.257859, -0.109325) (axis cs:0.247858, -0.107621) (axis cs:0.237861, -0.105902) (axis cs:0.227868, -0.104169) (axis cs:0.217879, -0.102420) (axis cs:0.207895, -0.100657) (axis cs:0.197914, -0.098878) (axis cs:0.187938, -0.097083) (axis cs:0.177966, -0.095271) (axis cs:0.167998, -0.093444) (axis cs:0.158035, -0.091599) (axis cs:0.148077, -0.089737) (axis cs:0.138123, -0.087857) (axis cs:0.128174, -0.085958) (axis cs:0.118231, -0.084041) (axis cs:0.108292, -0.082104) (axis cs:0.098359, -0.080148) (axis cs:0.088432, -0.078170) (axis cs:0.078510, -0.076172) (axis cs:0.068594, -0.074151) (axis cs:0.058684, -0.072108) (axis cs:0.048781, -0.070041) (axis cs:0.038883, -0.067951) (axis cs:0.028993, -0.065835) (axis cs:0.019109, -0.063694) (axis cs:0.009233, -0.061526) (axis cs:-0.000636, -0.059330) (axis cs:-0.010497, -0.057106) (axis cs:-0.020351, -0.054852) (axis cs:-0.030196, -0.052567) (axis cs:-0.040033, -0.050251) (axis cs:-0.049860, -0.047901) (axis cs:-0.059679, -0.045516) (axis cs:-0.069488, -0.043096) (axis cs:-0.079286, -0.040639) (axis cs:-0.089075, -0.038143) (axis cs:-0.098853, -0.035606) (axis cs:-0.108619, -0.033028) (axis cs:-0.118374, -0.030405) (axis cs:-0.128116, -0.027737) (axis cs:-0.137845, -0.025020) (axis cs:-0.147561, -0.022253) (axis cs:-0.157263, -0.019434) (axis cs:-0.166950, -0.016559) (axis cs:-0.176622, -0.013626) (axis cs:-0.186277, -0.010632) (axis cs:-0.195914, -0.007573) (axis cs:-0.205534, -0.004446) (axis cs:-0.215134, -0.001246) (axis cs:-0.224713, 0.002031) (axis cs:-0.234270, 0.005390) (axis cs:-0.243804, 0.008837) (axis cs:-0.253312, 0.012379) (axis cs:-0.262793, 0.016024) (axis cs:-0.272244, 0.019780) (axis cs:-0.281662, 0.023658) (axis cs:-0.291044, 0.027670) (axis cs:-0.300386, 0.031833) (axis cs:-0.309683, 0.036165) (axis cs:-0.318928, 0.040691) (axis cs:-0.328112, 0.045442) (axis cs:-0.337224, 0.050464) (axis cs:-0.346246, 0.055818) (axis cs:-0.355154, 0.061600) (axis cs:-0.363904, 0.067971) (axis cs:-0.372413, 0.075244) (axis cs:-0.380465, 0.084223) (axis cs:-0.385791, 0.103372)
};

\draw [smooth, ultra thin] plot coordinates {(axis cs:-0.050000, -0.325000) (axis cs:0.050000, -0.325000)
};

\draw [smooth, ultra thin] plot coordinates {(axis cs:0.000000, -0.580000) (axis cs:0.000000, -0.325000)
};

\draw [smooth, ultra thin] plot coordinates {(axis cs:-0.060000, -0.250000) (axis cs:-0.060000, -0.400000)
};

\draw [smooth, ultra thin] plot coordinates {(axis cs:-0.060000, -0.250000) (axis cs:0.060000, -0.250000)
};

\draw [smooth, ultra thin] plot coordinates {(axis cs:0.060000, -0.250000) (axis cs:0.060000, -0.400000)
};

\draw [smooth, ultra thin] plot coordinates {(axis cs:0.000000, -0.250000) (axis cs:0.000000, 0.000000)
};

\node[]    at    (axis cs:0.1, -0.025) {$m_s$};
\node[]    at    (axis cs:0.135, -0.325) {$c_s$};
\draw [smooth, solid, gray, ultra thin] plot coordinates {(axis cs:0.000600, 0.000000) (axis cs:-0.799400, 0.000000)
};

\draw [smooth, solid, gray, ultra thin] plot coordinates {(axis cs:0.000580, -0.000155) (axis cs:-0.772161, 0.206900)
};

\node[]    at    (axis cs:-0.6, 0.075) {$\theta(t)$};
\node[]    at    (axis cs:-0.58, -0.325) {$u_s$};
\draw [smooth, solid, gray, <->, ultra thin] plot coordinates {(axis cs:-0.650000, -0.580000) (axis cs:-0.650000, 0.000000)
};

\end{axis}
\end{tikzpicture}
 \caption{Foil-damper system} \label{FIG: FOIL DAMPER}
\end{figure}

In this section, the high-order, partitioned solver with the optimization framework introduced in
this document is applied to find the maximum energy harvesting through flow-induced oscillations of a 
NACA 0012 foil of length $l = 1$. The two-dimensional energy-harvesting model problem \cite{peng2009energy}
is represented by using a two-field FSI formulation. Consider the mass-damper system
in Figure~\ref{FIG: FOIL DAMPER}, the airfoil is suspended in an isentropic, viscous
flow where the rotational motion is a prescribed periodic motion
and the vertical displacement $u_s$ is determined by balancing the forces exerted
on the airfoil by the fluid and the damper (see \eqnref{EQ: SIMP STRUCT}).
The airfoil is initially at $ \theta(0) = 0 $, it matches a prescribed motion
for half a period, and then follows a periodic motion, as follows,
\begin{equation}
  \theta(t)=\begin{cases}
    \mu_A\cos(\frac{2t}{T}(\pi+ \mu_{\phi})), & t < \frac{T}{2}.\\
    \mu_A\cos(2\pi f t + \mu_{\phi}), & t \geq \frac{T}{2}.
  \end{cases}
\end{equation}
Here the period $T=5$ and the frequency is $f=0.2$.

The fluid is a perfect gas, with the adiabatic gas constant $\gamma = 1.4$, 
governed by the isentropic Navier-Stokes equations. The isentropic assumption states
the entropy of the system is assumed constant, which is tantamount to the flow
being adiabatic and reversible. For a perfect gas, the entropy is defined as
\begin{equation}\label{eqn:entropy}
  s = p/\rho^\gamma.
\end{equation}
The transformed conservation law, as described in Section~\ref{SUBSUBSEC: ALE Fluid}, 
is discretized with a standard high-order discontinuous
Galerkin method using Roe's flux \cite{roe1981approximate} for the inviscid
numerical flux and the Compact DG flux \cite{peraire2008compact} for the
viscous numerical flux. The DG discretization uses a mesh consisting of $3912$
cubic simplex elements ($p=3$). The second-order ODE in \eqnref{EQ: SIMP STRUCT} is
the governing equation for the mass-damper system with mass $m_s = 1$, damping
constant $c_s = 1$, stiffness $k_s = 0$, and external force given from the
fluid as described in Section~\ref{SUBSUBSEC: APP FSI STRUCT}. The mesh motion is
determined from the position and velocity of the structure using the
blending maps introduced in \cite{persson2009discontinuous} and identical to
that used in Section 5.1 of \cite{zahr2016adjoint}. IMEX4 is applied for temporal discretization, 
which matches the expected spatial order of accuracy obtained with polynomials
of degree 3.

The objective is to maximize the energy extraction $\mathcal{J} = \frac{1}{T}\int_T^{2T} c_s\dot{u_s}^2 dt $ by the device for the second period.
The energy injection to maintain the oscillation is defined by $E_{\theta} = -\frac{1}{T}\int_{T}^{2T} M_z\dot\theta dt$, where $M_z$ is the moment the fluid imparts 
onto the foil and $\dot\theta$ is the rotational speed of the foil.
We have linear constraints $-55^{\circ} \leq \mu_A^{\textrm{init}} \leq 55^{\circ} $ for the amplitude parameter and $-\frac{\pi}{2} < \mu_{\phi} < \frac{\pi}{2}$  for the phase parameter $\mu_{\phi}$, 
and a nonlinear constraint $E_{\theta} \geq -0.15$  for the energy injection $E_{\theta}$.

\begin{figure}
 \includegraphics[width=0.245\textwidth]{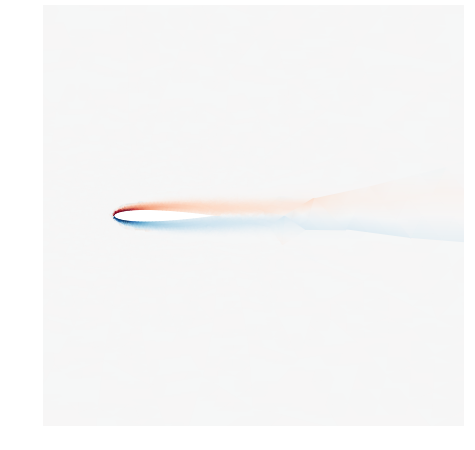}   \hspace{-0.1cm}
 \includegraphics[width=0.245\textwidth]{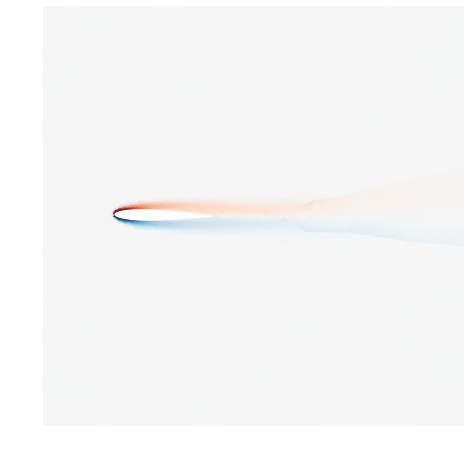}   \hspace{-0.1cm}
  \includegraphics[width=0.245\textwidth]{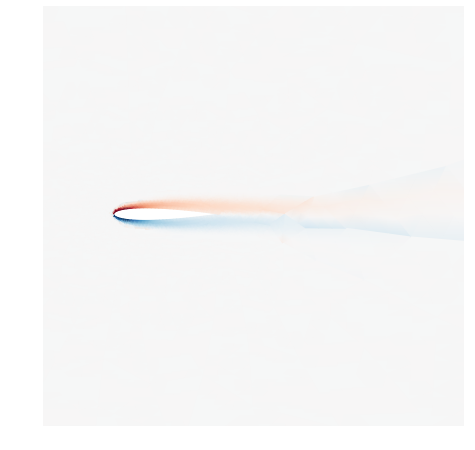}  \hspace{-0.1cm}
 \includegraphics[width=0.245\textwidth]{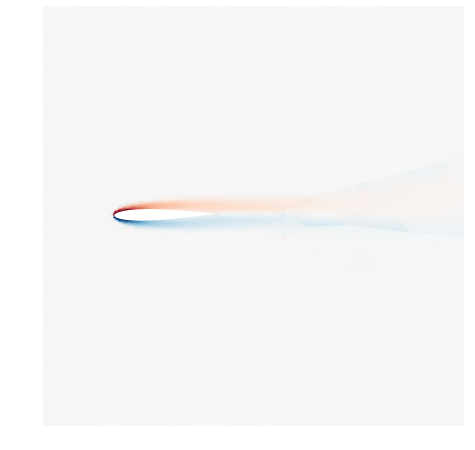} \\
 \caption{Airfoil motion and flow vorticity corresponding to foil-damper
          system under prescribed rotational motion
          $\theta(t) = \mu_A^{\textrm{init}}\cos(2\pi f t + \mu_{\phi}^{\textrm{init}})$ with frequency $f = 0.2$
          at various snapshots in time:
          $t = T,\,\frac{5}{4}T,\,\frac{6}{4}T,\,\frac{7}{4}T$
          (\emph{left}-to-\emph{right}, \emph{top}-to-\emph{bottom}).}
 \label{FIG: FOIL DAMPER VORT INIT}
\end{figure}

The initial motion is defined by  $\mu_A^{\textrm{init}} = 1^{\circ}$ and  $\mu_{\phi}^{\textrm{init}} = 0$. Snapshots of the vorticity
field and the motion of the airfoil are shown in Figure~\ref{FIG: FOIL DAMPER VORT INIT}, and the corresponding energy extraction is close to 0.

\begin{figure}
 \includegraphics[width=0.245\textwidth]{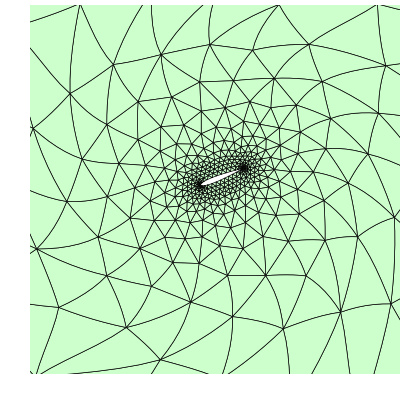}   \hspace{-0.1cm}
 \includegraphics[width=0.245\textwidth]{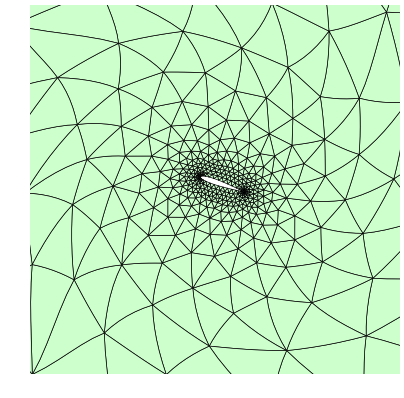}   \hspace{-0.1cm}
  \includegraphics[width=0.245\textwidth]{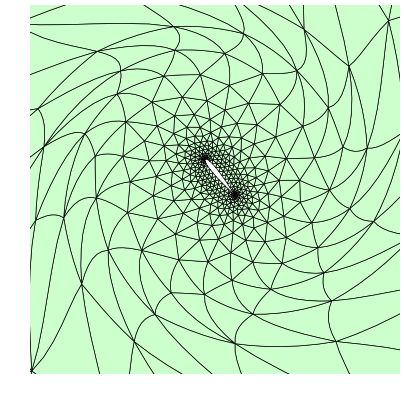}  \hspace{-0.1cm}
 \includegraphics[width=0.245\textwidth]{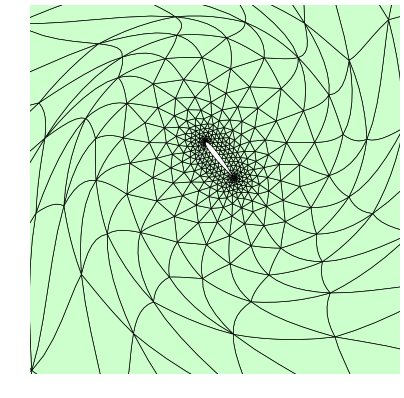} \\
 \includegraphics[width=0.245\textwidth]{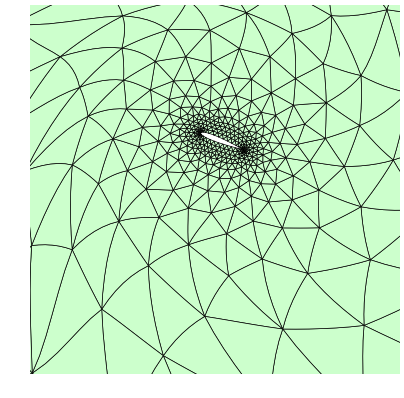}  \hspace{-0.1cm}
 \includegraphics[width=0.245\textwidth]{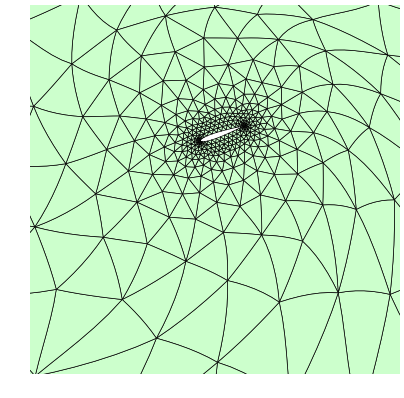}  \hspace{-0.1cm}
  \includegraphics[width=0.245\textwidth]{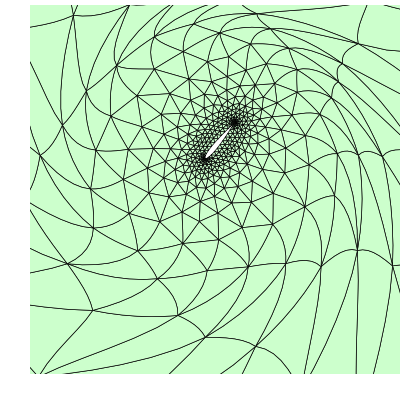}  \hspace{-0.1cm}
 \includegraphics[width=0.245\textwidth]{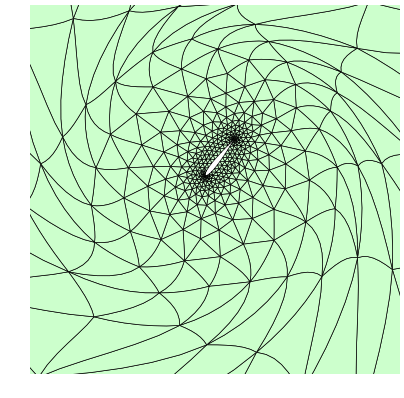}
 \caption{Airfoil motion and mesh deformation corresponding to foil-damper
          system under prescribed rotational motion
          $\theta(t) = \mu_A^{\textrm{opt}}\cos(2\pi f t + \mu_{\phi}^{\textrm{opt}})$ with frequency $f = 0.2$
          at various snapshots in time:
          $t = T,\,\frac{9}{8}T,\,\frac{10}{8}T,\,\frac{11}{8}T,\,\frac{12}{8}T,\,\frac{13}{8}T,\,\frac{14}{8}T,\,\frac{15}{8}T$
          (\emph{left}-to-\emph{right}, \emph{top}-to-\emph{bottom}).}
 \label{FIG: FOIL DAMPER MESH OPT}
\end{figure}

\begin{figure}
 \includegraphics[width=0.245\textwidth]{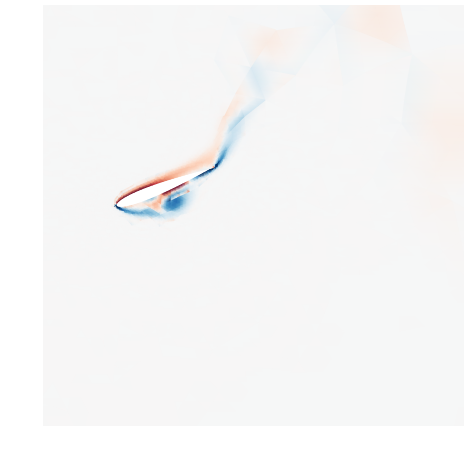}   \hspace{-0.1cm}
 \includegraphics[width=0.245\textwidth]{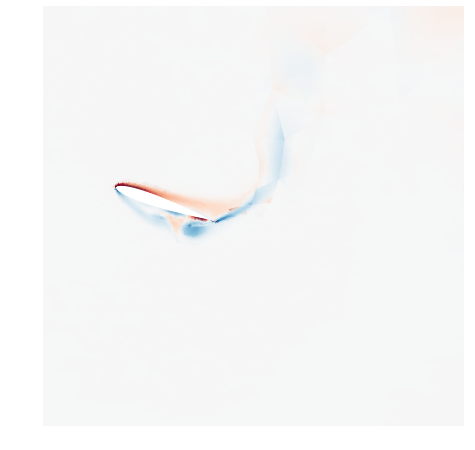}   \hspace{-0.1cm}
  \includegraphics[width=0.245\textwidth]{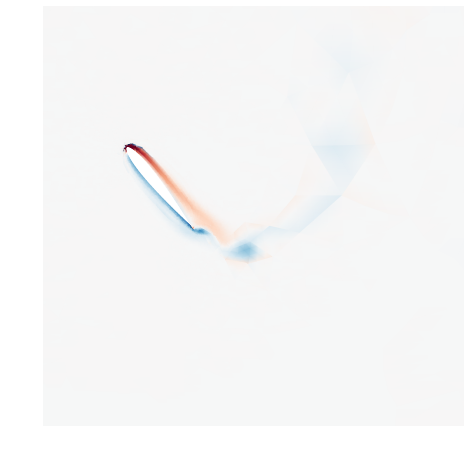}  \hspace{-0.1cm}
 \includegraphics[width=0.245\textwidth]{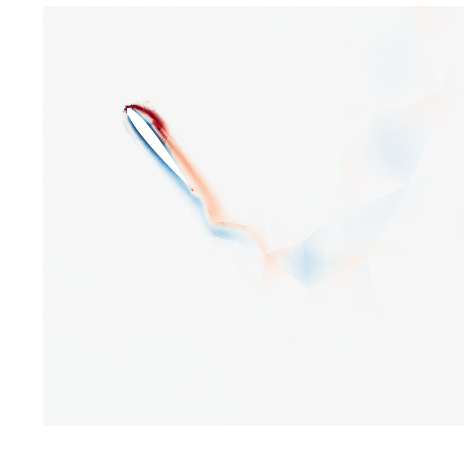} \\
 \includegraphics[width=0.245\textwidth]{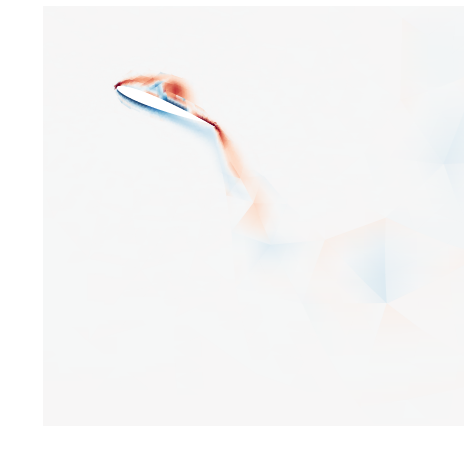}  \hspace{-0.1cm}
 \includegraphics[width=0.245\textwidth]{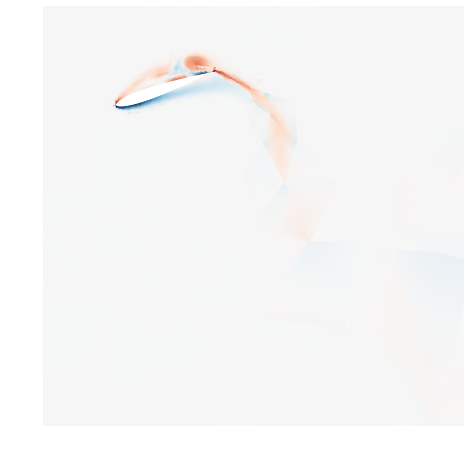}  \hspace{-0.1cm}
  \includegraphics[width=0.245\textwidth]{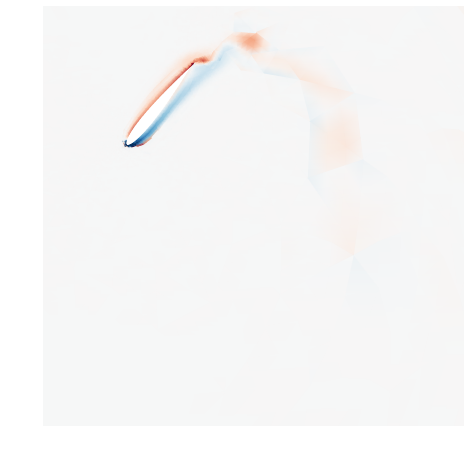}  \hspace{-0.1cm}
 \includegraphics[width=0.245\textwidth]{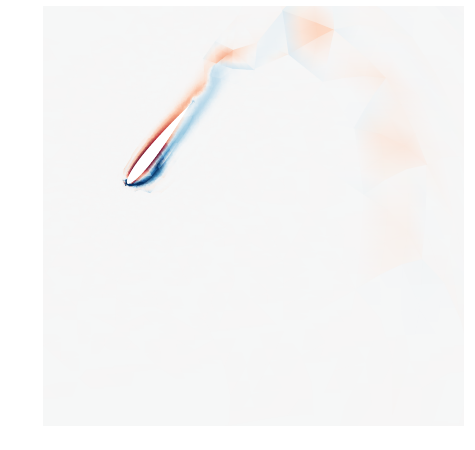}
 \caption{Airfoil motion and flow vorticity corresponding to foil-damper
          system under prescribed rotational motion
          $\theta(t) = \mu_A^{\textrm{opt}}\cos(2\pi f t + \mu_{\phi}^{\textrm{opt}})$ with frequency $f = 0.2$
          at various snapshots in time:
          $t = T,\,\frac{9}{8}T,\,\frac{10}{8}T,\,\frac{11}{8}T,\,\frac{12}{8}T,\,\frac{13}{8}T,\,\frac{14}{8}T,\,\frac{15}{8}T$
          (\emph{left}-to-\emph{right}, \emph{top}-to-\emph{bottom}).}
 \label{FIG: FOIL DAMPER VORT OPT}
\end{figure}

For the optimal oscillatory trajectory, the parameters obtained are  $\mu_A^{\textrm{opt}} = 55^{\circ}$ and  $\mu_{\phi}^{\textrm{opt}} = -22.95^{\circ}$.  
Snapshots of the ensued mesh motion are depicted in Figure~\ref{FIG: FOIL DAMPER MESH OPT}, snapshots of the vorticity
field and motion of the airfoil are shown in Figure~\ref{FIG: FOIL DAMPER VORT OPT}, and the energy extraction  $\mathcal{J} + E_{\theta}$
by this motion is almost 0.2.

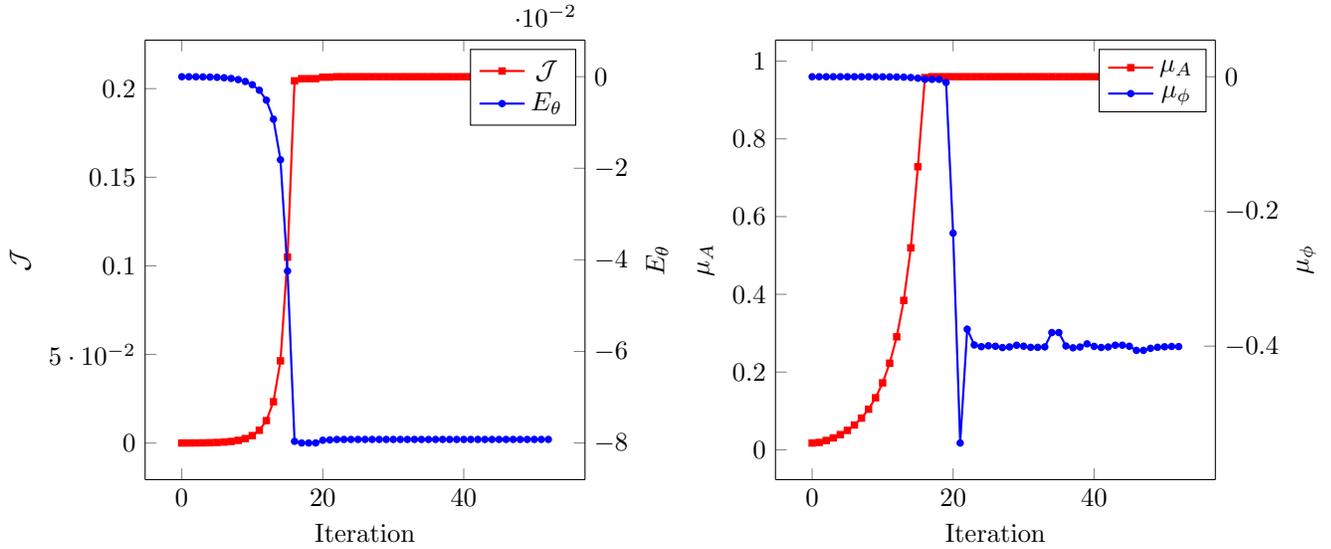
\begin{figure}
\begin{tikzpicture}
\begin{axis}[
    axis y line*=left,
    width=0.45\textwidth,
    height=0.45\textwidth,
    ylabel={$\mathcal{J}$},
    xlabel={Iteration}]
\addplot [red, solid, thick, mark=square*, mark size=1, mark options={solid}]  table[x index=0, y index=1] {data/naca_harvesting/naca_harvesting_optim_RK4.dat};\label{line:naca_harvesting_obj}
\end{axis}

\begin{axis}[
    axis y line*=right,
    axis x line=none,
    width=0.45\textwidth,
    height=0.45\textwidth,
    ylabel={$E_{\theta}$},
    xlabel={Iteration}]
\addlegendimage{/pgfplots/refstyle=line:naca_harvesting_obj}\addlegendentry{$\mathcal{J}$}
\addplot [blue, solid, thick, mark=otimes*, mark size=1, mark options={solid}]  table[x index=0, y index=2] {data/naca_harvesting/naca_harvesting_optim_RK4.dat};\label{line:naca_harvesting_constraint}
\addlegendentry{$E_{\theta}$}
\end{axis}

\end{tikzpicture}
\begin{tikzpicture}
\begin{axis}[
    axis y line*=left,
    width=0.45\textwidth,
    height=0.45\textwidth,
    ylabel={$\mu_{A}$},
    xlabel={Iteration}]
\addplot [red, solid, thick, mark=square*, mark size=1, mark options={solid}]  table[x index=0, y index=3] {data/naca_harvesting/naca_harvesting_optim_RK4.dat};\label{line:naca_harvesting_mu_A}
\end{axis}

\begin{axis}[
    axis y line*=right,
    axis x line=none,
    width=0.45\textwidth,
    height=0.45\textwidth,
    ylabel={$\mu_{\phi}$},
    xlabel={Iteration}]
\addlegendimage{/pgfplots/refstyle=line:naca_harvesting_mu_A}\addlegendentry{$\mu_A$}
\addplot [blue, solid, thick, mark=otimes*, mark size=1, mark options={solid}]  table[x index=0, y index=4] {data/naca_harvesting/naca_harvesting_optim_RK4.dat};\label{line:naca_harvesting_mu_phi}
\addlegendentry{$\mu_{\phi}$}
\end{axis}

\end{tikzpicture}
\caption{Convergence of the optimizer for the NACA harvesting problem.}
\label{FIG: NACA HARVESTING}
\end{figure}

The convergence of the objective function $\mathcal{J}$ and the nonlinear constraint $E_{\theta}$ are reported in Figure~\ref{FIG: NACA HARVESTING}-left. 
The convergence of the parameters $\mu_A$ and $\mu_{\phi}$ are presented in Figure~\ref{FIG: NACA HARVESTING}-right. Initially, the energy harvester extract almost no energy from the fluid without energy injection. However, for the optimal oscillatory trajectory $\theta(t) = \mu_A\cos(2\pi f t + \mu_{\phi})$, the injected energy for maintaining the oscillation is $E_{\theta} = -7.92\times 10^{-2}$; The energy extracted by the damper is $\mathcal{J} = 2.07\times 10^{-1}$.  The optimized energy harvester can extract $\mathcal{J} + E_{\theta} = 1.27\times 10^{-1}$ from the fluid flow, which demonstrate the potential benefits of multiphysics optimization.

\section{Conclusion}\label{SEC: CONCLUSION}

We have presented a framework for optimizing unsteady multiphysics systems, based on the high-order, linearly stable, 
partitioned solver introduced in \cite{huang2018high}. 
An implicit-explicit Runge-Kutta scheme was used for high-order temporal
integration with the benefit of achieving accuracy beyond second-order and decoupling all subsystems.
Therefore, the corresponding adjoint equations or sensitivity equations can be solved in a partitioned 
manner, i.e. subsystem-by-subsystem and substage-by-substage. 
While we did not quantify the benefits of high-order discretizations for these optimization problems, it is still likely that high-order spatial and temporal accuracy 
allow for smaller mesh size and larger timestep size, which improve the efficiency of function and gradient evaluations in the optimization procedure. 
Due to the fully discrete adjoint solver, exact gradients are obtained, and the implementation was verified using finite differences. A gradient-based optimizer converged quickly to optional solutions for our examples problems.
In future work, the efficiency of the present optimization framework will be studied, and it will be used to better understand the energy harvesting process with multiple airfoils and for the optimization of 3D fluid-structure systems.  

\section*{Acknowledgements}

This work was supported in part by the Luis W. Alvarez Postdoctoral Fellowship (MZ), by the Director,
Office of Science, Office of Advanced Scientific Computing Research, U.S. Department of Energy under
Contract No. DE-AC02-05CH11231 (MZ, PP),  by the NASA National Aeronautics and Space Administration
under grant number NNX16AP15A (MZ, PP),  by the Jet Propulsion Laboratory (JPL) under Contract JPL-RSA No. 1590208  (DH)
, and by the National Aeronautics and Space Administration
(NASA) under Early Stage Innovations (ESI) Grant NASA-NNX17AD02G (DH). The content of this publication does not necessarily
reflect the position or policy of any of these supporters, and no official endorsement should be inferred.

\bibliographystyle{unsrt}
\bibliography{myref}
\end{document}